\documentclass[12pt]{article}

\def\be{\begin{equation}}
\def\ee{\end{equation}}
\def\bea{\begin{eqnarray}}
\def\eea{\end{eqnarray}}
\def\bean{\begin{eqnarray*}}
\def\eean{\end{eqnarray*}}

\def\S{\Sigma}
\def\P{\Psi}
\def\G{\Gamma}
\def\D{\Delta}
\def\U{\Upsilon}

\def\square{\hfill\hbox{\vrule\vbox{\hrule\phantom{N}\hrule}\vrule}\,}

\newtheorem{defi}{Definition}[section]

\newtheorem{theorem}{Theorem}[section]
\newtheorem{result}{Result}[section]

\newtheorem{coro}{Corollary}[section]

\begin{document}

\title{A weighted de Rham operator acting on arbitrary tensor fields and their local potentials.}
\author{S. Brian Edgar$^1$ and Jos\'e M.M. Senovilla$^2$
\\ $^1$ Matematiska institutionen,   Link\"{o}pings universitet\\ 
Link\"{o}ping, Sweden S-581 83\\
$^2$ F\'{\i}sica Te\'orica, Universidad del Pa\'{\i}s Vasco, \\
Apartado 644, 48080 Bilbao, Spain\\ 
e-mails: bredg@mai.liu.se, josemm.senovilla@ehu.es}

%\date{}
\maketitle

\begin{abstract} 
We introduce a {\it weighted de Rham operator} which acts on 
arbitrary tensor fields by considering their structure as $r$-fold forms.
We can thereby define {\it associated superpotentials} for all tensor 
fields in all dimensions and, from any of these superpotentials, we
deduce in a straightforward and natural manner the existence of $2r$
potentials for any tensor field, where $r$ is its form-structure
number.  By specialising this result to {\em symmetric} double forms, we
are able to obtain a pair of potentials for the Riemann tensor, and a
single $(2,3)$-form potential for the Weyl tensor due to its
tracelessness.  This latter
potential is the $n$-dimensional version of the double dual of the
classical four dimensional $(2,1)$-form Lanczos potential.  We also
introduce a new concept of {\em harmonic} tensor fields,
demonstrate that the new weighted de Rham operator has many other 
desirable properties and, in particular, it is the natural
operator to use in the Laplace-like equation for the Riemann tensor.
\end{abstract}


\


\


\section{Introduction}

The classical Helmholtz theorem \cite{H} characterises any three
dimensional vector field $\vec V$ on Euclidean space
in terms of two potentials (scalar
and vector respectively), and the proof depends on the existence of a
solution to the vector version of Poisson's equation
$$\nabla^2 \stackrel{o}{\vec V} = {\vec V}$$ 
for a superpotential vector field $\stackrel{o}{\vec V}$. 
The potentials follow immediately from the vector operator identity
$${\vec V}=\nabla^2{\stackrel{o}{\vec V}} \equiv \nabla 
({\nabla }\cdot \stackrel{o}{\vec V}) - 
\nabla \times( \nabla \times \stackrel{o}{\vec V}).$$
Four dimensional \cite{W} and $n$-dimensional generalisations
\cite{LD},\cite{Gem} have been discussed for vector fields using
analogous arguments.

More generally, on an $n$-dimensional  pseudo-Riemannian manifold, 
the (local) Hodge decomposition \cite{Ho} (also \cite{Fl,MTW,Na,GS,Fr,Sc}) 
characterises any $p$-form $\Sigma$  
in terms of a pair of potentials (respectively a $(p-1)$-form and a 
$(p+1)$-form), and the proof depends on the existence of a solution to the 
Laplace-like equation
$$\Delta\stackrel{o}{\mbox{$\Sigma$}} = \Sigma$$ 
for a superpotential $\stackrel{o}{\mbox{$\Sigma$}}$, 
where $\Delta \equiv d\delta + \delta d$ is the de Rham operator for 
$p$-forms \cite{DR}, $d$ is the  exterior differential  
and $\delta$ is the codifferential operator.
The potentials follow immediately (see Result \ref{r2}) from this definition 
$$\S= \Delta \stackrel{o}{\mbox{$\Sigma$}}\ \equiv 
d (\delta \stackrel{o}{\mbox{$\Sigma$}})+
\delta (d \stackrel{o}{\mbox{$\Sigma$}}) \ .$$

In the case of general tensor fields on an $n$-dimensional 
pseudo-Riemannian manifold with Levi-Civita connection $\nabla$, one
is still able to exploit the natural generalisations of these three
operators  by considering the tensor as an appropriate {\it
tensor-valued $p$-form}, e.g. \cite{Bel,MTW,Coq}.  
It is important to note that when a tensor (with more than $p$ 
antisymmetrical indices) is considered
as a tensor-valued $p$-form, the generalised de Rham operator 
 acting on such a tensor does not
commute with properties involving indices from {\it both} its 
$p$-form part and its tensor-valued part; hence simple natural constructions
for a superpotential and potentials for a tensor-valued $p$-form
cannot be deduced in an analogous manner as for single $p$-forms.

Recognising this deficiency,
Lichnerowicz \cite{Lic0,Lic} proposed yet another generalisation $\Delta_L$ which,
when acting on {\it arbitrary tensors}, commutes with all their index
properties; unfortunately, the Lichnerowicz operator $\Delta_L$ does
not have direct links to the exterior differential $d$ and
codifferential operator $\delta$, and hence to natural potentials in
the manner of $\D=d\delta+\delta d$ for single $p$-forms.

In this paper we overcome all these deficiencies in a different manner 
by exploiting the fact that {\em any} tensor
can be considered, in a precise and unequivocal way,
as an {\it $r$-fold form} \cite{S}, and defining a
new {generalised Laplacian} which is in effect a {\it weighted de Rham
operator} $\bar\Delta$; this operator is directly suited to the $r$-fold
form structure of the given tensor and commutes with all index properties, 
but crucially also has
direct links with appropriate exterior differential $d_{(i)}$ and codifferential
$\delta_{(i)}$ operators (see Section \ref{sec:r-fold}).  Using the new operator 
$\bar \Delta$ we are able to
identify associated superpotentials for any tensor directly via a
simple Laplace-like equation.  Furthermore, we obtain, in an analogous
manner, the natural generalisation of the Helmholtz and Hodge
decomposition: in arbitrary curved spaces, when an arbitrary tensor is
considered as an $r$-fold form then it can be expressed in a very
natural and useful manner in terms of $2r$ local potentials.  By
specialising these results to {\it double forms}, we can then deduce
the potential structure of Riemann and Weyl curvature tensors
respectively, and hence re-establish our recent result which obtained
one double $(2,3)$-form potential for the Weyl tensor in all
dimensions \cite{ES}.

We begin in Section 2 by formulating familiar results for $p$-forms in
a manner from which we can generalise to more general tensors, and in
Section 3Ê we first review how to consider any tensor $T$ as an
$r$-fold form \cite{S}.  Associated with each block, $i=1,2, \ldots r$
of an $r$-fold form $T$ we have the three generalised operators $d_{(i)}, \delta_{(i)},
\D_{(i)}$ from which we are able to determine a superpotential and a pair
of potentials; but these superpotentials and potentials are highly
non-unique with limited practical use.  With the introduction of the
{\it weighted de Rham operator} 
$$\bar\D \equiv \frac{1}{r}(\D_{(1)} + \D_{(2)} +
\ldots \D_r)=\frac{1}{r}\sum_{i=1}^r \D_{(i)}$$ 
in Section 4, we establish our basic result that for an
$r$-fold form $T$ there exists an {\it associated local superpotential}
$ {\stackrel{o}{\mbox{$T$}}}$ given by 
$$\bar \Delta \stackrel{o}{\mbox{$T$}}=T,
$$
where it is important to note that this superpotential $
{\stackrel{o}{\mbox{$T$}}}$ does not just have the same form structure, but
also has the same index symmetry and trace properties as $T$.  From
this we can deduce, in a natural and straightforward manner, a
generalisation of the Hodge decomposition, demonstrating the existence
of $2r$ local potentials, 
$$T =\bar\D \stackrel{o}{\mbox{$T$}} \  =
\frac{1}{r}\sum_{i=1}^r \D_{(i)}\stackrel{o}{\mbox{$T$}}\   \equiv
\frac{1}{r}\sum_{i=1}^r
[\delta_{(i)}(d_{(i)}\stackrel{o}{\mbox{$T$}})+
d_{(i)}(\delta_{(i)}\stackrel{o}{\mbox{$T$}})] \, .
$$
This leads naturally to  a new definition for harmonic tensors, in general.

In Section 5 we specialise these results to double $(q,p)$-forms and
identify their four different potentials.  Further specialisations
are made in Section 6: first to (anti)symmetric double $(p,p)$-forms, 
which are shown to have only a pair of potentials --- a double 
$(p,p+1)$-form ÊÊ
and a double $(p,p-1)$-form --- and then to the particular 
value $p=2$. These results are immediately 
applicable to the Riemann curvature tensor---Theorem \ref{theoremR}--- and to 
the traceless Weyl tensor ---Theorem \ref{theoremW}, \cite{ES}.

In Section 7 we demonstrate that the new weighted de Rham operator is
the natural operator to use in the Laplace-like equation for the
Riemann tensor, and we illustrate its advantages over the particular
de Rham operator which is usually used.  In Section 8 we summarise, demonstrate how this work  provides insights into some other investigations, and
outline future plans.

\section{Standard results for $p$-forms.}
Let $V_n$ be any differentiable $n$-dimensional pseudo-Riemannian
manifold endowed with a metric $g_{ab}$ of arbitrary signature.  The
covariant derivative associated with $g$ is denoted either by $\nabla$
or by a semicolon, e.g., $\nabla_{a}v_{b}\equiv v_{b;a}$; and its
Riemann and Ricci tensors by $R_{abcd}$ and $R_{ab}\equiv
R^{c}{}_{acb}$, respectively.  Our convention for the Riemann tensor
follows from the Ricci identity:
$$
2v_{a;[bc]}=R^d{}_{abc}v_{d} .
$$
As usual, we use round and square brackets 
to indicate symmetrization and antisymmetrization of indices, respectively.

The graded algebra of exterior forms is $\Lambda$, with $\Lambda^p$ denoting the 
set of exterior $p$-forms; in particular, for $p=0$, the elements of $\Lambda^0$ will be scalar functions, (also called sometimes $0$-forms).  The canonical volume element $n$-form 
$\eta\in \Lambda^n$ is denoted by 
$\eta_{a_{1}\dots a_n}=\eta_{[a_1\dots a_n]}$. Then, we can define the 
standard Hodge dual operator $* : \Lambda^p \longrightarrow \Lambda^{n-p}$ by 
means of
\be
\stackrel{*}{\S}_{a_{p+1}\dots a_{n}}\equiv
\frac{1}{p!}\eta_{a_1\dots a_n}\S^{a_1\dots a_p} \hspace{1cm} \forall\, \S\in\Lambda^p .
\label{*}
\ee
It is easy to prove that $**\,=\epsilon (-1)^{p(n-p)}$ where 
$\epsilon =\pm 1=$sign(det($g_{ab}$)) is a sign depending on the signature.

One can define a scalar product $<\, ,\, >$ on each $\Lambda^p$ in the standard way by 
first defining the function
$$
(\S,\Phi)\equiv \S_{a_{1}\dots a_{p}}\Phi^{a_{1}\dots a_{p}}
$$
for every $\S,\Phi\in\Lambda^p$ and then integrating this over the 
manifold $V_{n}$ 
\be
<\S ,\Phi >\,  \equiv \int_{V_{n}} (\S,\Phi)\, \eta = p!
\int_{V_{n}} \S\, \wedge \stackrel{*}{\mbox{$\Phi$}} \, .\label{scalar-product}
\ee
When $V_{n}$ is compact without boundary the integration can be
carried out without further restrictions; otherwise, at least one of
$\S,\Phi$ has to be assumed to have compact support.  This will be
assumed without explicit mention in what follows.  This scalar product
is bi-linear, symmetric and non-degenerate.  In the case of proper
Riemannian manifolds, it is also positive definite.

The standard operations on the exterior algebra are the exterior 
differential (also called `curl') $d\, : \Lambda^p \longrightarrow \Lambda^{p+1}$ and 
codifferential (also called 'divergence')
$\delta \, : \Lambda^p \longrightarrow \Lambda^{p-1}$ \cite{Car, Lic, Fl, BT,MTW,Th}.
The second is simply defined by 
\be
\delta \equiv (-1)^p *^{-1}d\, * = \epsilon\, (-1)^{(n-p)(p-1)+1} * d\, *
\label{delta}
\ee
when acting on $p$-forms. With index notation they can be given as
\bean
(d\S)_{a_{1}\dots a_{p+1}}\equiv (p+1)\nabla_{[a_{1}}\S_{a_{2}\dots 
a_{p+1}]}=(-1)^p (p+1)\S_{[a_{1}\dots a_{p};a_{p+1}]} \, ,\\
(\delta\S)_{a_{2}\dots a_{p}}\equiv - \nabla^{a_{1}}\, \S_{a_{1}a_{2}\dots a_{p}}=
(-1)^{p}\S_{a_{2}\dots a_{p}a_{1}}{}^{;a_{1}} \hspace{1cm}
\eean
for all $\S\in\Lambda^p$. Observe that $d \alpha =0,\,\,\, \forall 
\alpha\in\Lambda^n$ and, as  is implicit in this formula 
and follows trivially from the definition (\ref{delta}) of $\delta$, for 
functions $f$  one has: $\delta 
f = 0,\,\,\, \forall f\in \Lambda^0$.
Elementary properties of these operators are $d^2 \equiv 0$ and 
$\delta^2 \equiv 0$. Moreover, for every $\S\in\Lambda^p$ and 
$\G\in\Lambda^{p+1}$ we have
$$
<d\S ,\G >\, =\, <\S , \delta \G >
$$
as can be easily checked by using the Gauss theorem. Therefore, $d$ 
and $\delta$ are mutually adjoint with respect to $<\, ,\, >$ 
\cite{Lic0,Bel}. As usual, this property can be used to define $d$ 
and $\delta$ acting on $p$-form distributions by means of the standard 
use of ``test $p$-forms'' (i.e., $C^{\infty}$ $p$-forms of 
compact support), see \cite{Lic0}.

A $p$-form $\S$ is called {\em closed} if $d\S=0$ 
and {\em exact} if $\S=d \P$ for some $\P\in\Lambda^{p-1}$. Analogously $\S$ is 
said to be {\em co-closed} if $\delta\S=0$ and {\em co-exact} if $\S=\delta \G$ 
for some $\G\in \Lambda^{p+1}$. In general, every statement on $p$-forms 
has a dual statement replacing $d$ for $\delta$ and the form by its 
Hodge dual. Obviously, all exact $p$-forms are closed, but the 
converse is not true, globally, in general. The Poincar\'e 
lemma ensures that every closed form is locally exact. Analogously, 
all co-exact forms are co-closed, and these are locally co-exact. 

The de Rham cohomology class of order $p$ is defined as the quotient of 
the set of closed $p$-forms by the set of exact $p$-forms. And the de 
Rham Laplacian operator $\D :\Lambda^p \longrightarrow \Lambda^p$ is 
defined intrinsically by $\D \equiv d\delta + \delta d$ 
(sometimes also written as $\D = (d +\delta)^2$) \cite{DR}. Its formula with 
index notation is
\be
(\D \S)_{a_{1}\dots a_{p}}=-\nabla^c\nabla_{c}\S_{a_{1}\dots 
a_{p}}+pR_{c[a_{1}}\S^c{}_{a_{2}\dots 
a_{p}]}-\frac{p(p-1)}{2}R_{cd[a_{1}a_{2}}\S^{cd}{}_{a_{3}\dots 
a_{p}]}\, .
\label{deR}
\ee
Observe that the Riemann tensor terms 
are not present for 1-forms, while for 0-forms $\D$ is (minus) the usual 
Laplacian $\nabla^c\nabla_{c}=g^{bc}\nabla_{b}\nabla_{c}$.
%% 
 % Note that since $\Sigma$ and $\Delta \Sigma$ are both forms of the
 % same type, and that the operator $\Delta$ commutes with traces, we can
 % remove the brackets and use the notation $ \D \S_{a_{1}\dots a_{p}}$
 % for the index form of the de Rham Laplacian without any danger of
 % ambiguity.
 %%

The following important properties of $\D$ are easily obtained (see 
e.g. \cite{DR,Lic0,Fl,Na}) 
\begin{result}
    \label{r1}
The operator $\D$ is linear, self-adjoint with respect to $<\, ,\, >$
\be
<\D \S , \Phi >\, =\, <\S , \D \Phi > \hspace{3mm} \forall \S,\Phi \in \Lambda^p
\ee
and also commutes with $*$, $d$ and $\delta$
\be
*\D = \D * , \hspace{1cm}d\Delta=\Delta d , \hspace{1cm} 
\delta\Delta = \Delta\delta \, .
\label{comm}
\ee
\end{result}
Furthermore, we can immediately  prove the identity
\be
<\S ,\D\S >\, =\, <d\S ,d\S > + <\delta\S ,\delta\S> \, .\label{dotp}
\ee
In the case of proper Riemannian manifolds this implies that 
$\D$ is a positive operator, i.e., $<\S ,\D\S >\,\,  \geq 0$ for all 
$\S\in\Lambda^p$.

In what follows, and in order to keep the technical complications to a 
minimum, we shall assume that the metric $g$ is analytical. This will 
allow us to use the simpler versions of the existence results for 
solutions of Laplace-like equations. More powerful results could of 
course be used by assuming $C^2$ differentiability of the metric, but 
this will not change the main points we wish to make in this paper. 
Keeping this in mind, a key result for $p$-forms \cite{Fl,Na,Fr}, 
which we will generalise in this paper, is
\begin{result}[Local Helmholtz-Hodge decomposition]
    \label{r2}
Given any $p$-form $\S \in\Lambda^{p} $, there always exists a local superpotential 
$\stackrel{o}{\mbox{$\Sigma$}} \  \in\Lambda^{p} $ such that
$\D  \stackrel{o}{\mbox{$\Sigma$}}\, = \S $; 
furthermore there always exists a pair of local 
potentials $(\P,\G)$ with $\P\in\Lambda^{p-1}$ and $\G\in\Lambda^{p+1}$ such 
that 
\be
\S=d\P + \delta\G \label{r2eq}
\ee
where $\P=\delta  \stackrel{o}{\mbox{$\Sigma$}}$ and 
$\G=d  \stackrel{o}{\mbox{$\Sigma$}}$.    
\end{result}

\begin{proof}  
From the structure of the de Rham operator $\D$ in (\ref{deR}),
according to the Cauchy-Kovalewski theorem \cite{CH}, we can always
find a local solution $\stackrel{o}{\Sigma}$, which is another $p$-form,
to the equation $\S=\D \stackrel{o}{\mbox{$\Sigma$}}$.  The
remainder of the result then follows from the definition 
$\D \equiv d\delta + \delta d$. \square
\end{proof}

Note that $\Psi$ is missing for the case $p=0$ so that in this case the 
statement is just that any function is locally the divergence of 
a vector field, while for $p=n$, $\Gamma$ is identically zero stating 
that any $n$-form is closed ergo locally exact.

An index version of Result \ref{r2} is
\bea
(\D{ \stackrel{o}{\mbox{$\Sigma$}}})_{a_{1}\dots a_{p}}& = &
(-1)^{p}p{ \stackrel{o}{\mbox{$\Sigma$}}}_{c[a_{1}\dots a_{p-1}}{}^{;c}{}_{a_{p}]}
-(p+1){ \stackrel{o}{\mbox{$\Sigma$}}}_{[a_{1}\dots a_{p};c]}{}^c \nonumber \\
&=&(-1)^{p-1} p \P_{[a_{1}\dots a_{p-1};a_{p}]}-
\G_{ca_{1}\dots a_{p}}{}^{;c}\nonumber \\ 
& =& (d\P)_{a_{1}\dots a_{p}}+
(\delta\G)_{a_{1}\dots a_{p}} \nonumber \\ 
&= &\S_{a_{1}\dots a_{p}} \label{calcul}
\eea
where we have defined 
\bea
\P_{a_{1}\dots a_{p-1}}&\equiv& - 
\stackrel{o}{\mbox{$\Sigma$}}_{ca_{1}\dots a_{p-1}}{}^{;c}
%\left[(-1)^{p-2}(p-1)\Phi_{[a_{1}\dots a_{p-2};a_{p-1}]}\right]
\  \  , \\
\G_{a_{1}\dots a_{p+1}}& \equiv& (-1)^p(p+1)  
\stackrel{o}{\mbox{$\Sigma$}}_{[a_{1}\dots a_{p};a_{p+1}]} \, .
%\left[\epsilon \Pi_{ca_{1}\dots a_{p+1}}{}^{;c}\right]
\eea
Observe that this result is independent of any field equations for $\S$, and 
is a {\em local} result.  Actually, Result \ref{r2} can be 
strengthened in the case of compact without 
boundary proper Riemannian manifolds and one can obtain 
the {\em global Hodge decomposition theorem} \cite{Ho,DR,Gol,Lic0,Fl,Na,Fr,Sc}: any 
$\S\in\Lambda^p$ admits a {\em unique} global decomposition as $\S =d\P 
+\delta\G +\U$ where $\P\in\Lambda^{p-1}$, $\G\in\Lambda^{p+1}$ and $\U\in\Lambda^p$ 
is a {\em harmonic} $p$-form (that is, $\U$ is closed and co-closed, 
i.e.,  $\D \U =0$) in the same co-homology class as $\S$. Note that the global 
Hodge result can be used to obtain the more modest local Result \ref{r2} because, 
since  $\U$ is closed and co-closed, it is also locally exact and 
co-exact (see \cite{Gol}).

Closely linked to the existence of potentials is their gauge freedom, 
and there are well known results for $p$-forms which exploit this gauge
freedom. However, when one considers
the more complicated {\it tensor-valued forms} which are the subject
of this paper, there are aspects of the role of gauge which are
significantly different from the $p$-form results.  
A detailed discussion of the role and application of
gauge for tensor valued forms will be presented in \cite{ESg}.

\section{Generalization to arbitrary tensors: $r$-fold forms}
\label{sec:r-fold}
In the previous section the two crucial results were being able to
identify a superpotential via a Laplace-like equation for the de Rham
operator, and then being able to link potentials in a simple natural
manner to derivatives of the superpotential, via this de Rham
operator, which is of course defined only for differential forms.  The
question now is how to generalize all this to arbitrary tensors.

Long ago Lichnerowicz \cite{Lic0,Lic} proposed a generalised Laplacian
for arbitrary tensor fields which had the first of these properties,
as well as a number of other useful properties.  Given an arbitrary
rank-$m$ tensor field $T_{a_{1}\dots a_{m}}$, the Lichnerowicz
operator $\D_L$ can be defined as
\bea
(\D_{L} T)_{a_{1}\dots a_{m}}\equiv - \nabla^c\nabla_{c} T_{a_{1}\dots 
a_{m}}&+& \sum_{s=1}^m R^c{}_{a_{s}}T_{a_{1}\dots a_{s-1}ca_{s+1}\dots 
a_{m}}\\ \nonumber
 & & - \sum_{s\neq t}^m R^c{}_{a_{s}}{}^d{}_{a_{t}}
T_{a_{1}\dots a_{s-1}ca_{s+1}\dots a_{t-1}da_{t+1}\dots a_{m}}\, .
\label{Li}
\eea
This operator has a number of very important properties \cite{Lic0,Lic}:
\begin{itemize}
\item $\Delta_L$ respects the symmetry 
properties of $T_{a_{1}\dots a_{m}}$; that is to say, 
$(\D_{L} T)_{a_{1}\dots a_{m}}$ has exactly the same index symmetries 
as $T_{a_{1}\dots a_{m}}$. 
\item $\D_{L}$ commutes with traces, i.e., 
 the trace on any two indices of $(\D_{L} T)_{a_{1}\dots 
a_{m}}$ equals 
$\D_{L}$ applied to the corresponding trace of $T_{a_{1}\dots a_{m}}$; 
in particular, if $T_{a_{1}\dots a_{m}}$ is 
traceless in any pair of indices, then so is $(\D_{L} T)_{a_{1}\dots 
a_{m}}$.
\item $\Delta_L$ is self-adjoint with respect to the scalar product 
$\{\, ,\, \}$ defined by
\be
\{T,S\}\, \equiv \int_{V_{n}} T_{a_{1}\dots a_{m}}S^{a_{1}\dots 
a_{m}}\, \eta \label{scalarL}
\ee
for arbitrary $T,S\in T_{m}(V_{n})$.  (Of course, the same comments as 
followed (\ref{scalar-product}) on 
the compactness of $V_{n}$ or of the support of one of the tensors, are 
in order here.) Therefore $\{\D_{L}T,S\}=\{T,\D_{L}S\}$.
\item In manifolds with a parallel Ricci tensor ($\nabla_{a}R_{bc}=0$, 
which includes the important cases of Einstein spaces and Ricci-flat manifolds,)
the following two properties hold; see \cite{Lic0,Lic}:
\begin{itemize}
\item when acting on rank-1 tensors, $\Delta_L$ commutes with the 
covariant derivative;
\item when acting on rank-2 tensors, $\D_L$ commutes with the divergence operator.
\end{itemize} 
\item {\em When acting on $p$-forms}, $\Delta_L$ coincides  with 
the de Rham Laplacian $\Delta$:
$$
\D_{L}\S = \D\S , \hspace{1cm} \forall \,\, \S\in \Lambda^p \, .
$$
\end{itemize}

Unfortunately, the Lichnerowicz Laplacian $\D_L$ does not have simple
natural links to first derivative operators for arbitrary tensors in
the same way as the de Rham operator has for $p$-forms.

We believe that the 
best way to extend the principles of Section 2 to arbitrary tensors 
is to consider tensors as 
{\em $r$-fold forms}. This terminology was extensively considered in
\cite{S}, and the underlying simple idea is based on the following 
remark: Given any rank-$m$ tensor $T^{a_{1}\dots a_{m}}$, there is a 
{\em minimum} natural number $r$, $r\leq m$ and a unique set of $r$ 
natural numbers $n_{1},\dots ,n_{r}$, with $\sum_{i=1}^{r}n_{i}=m$, 
such that $T^{a_{1}\dots a_{m}}$ is a linear 
map on $\Lambda^{n_{1}}\times\dots\times \Lambda^{n_{r}}$. In other words, 
there always exists a minimum $r$ such that
$\tilde{T}\in \Lambda_{n_{1}}\otimes\dots\otimes \Lambda_{n_{r}}$, 
where $\tilde{T}^{a_{1}\dots a_{m}}$ is the appropriate permuted version of 
$T^{a_{1}\dots a_{m}}$ which 
selects the natural order for the $n_{1},\dots ,n_{r}$ entries. Tensors 
seen in this way are called $r$-fold $(n_{1},\dots ,n_{r})$-forms 
\cite{S}. In short, {\em all tensors can be considered, in a precise 
way, as $r$-fold forms}. 
\begin{defi}[Form-structure number and block ranks]
For any tensor $T$, the uniquely defined number $r$ will be called 
its {\em form-structure number}, and each of the $n_{i}$, the {\em $i$-th block
rank}.
\end{defi}

Some simple examples are: any $p$-form $\S$ is trivially a 
single (that is, 1-fold) $p$-form, while $\nabla\S$ is a double 
$(1,p)$-form, with $r=2$ and $n_{1}=1, n_{2}=p$; 
the Riemann tensor ($r=2$, $n_{1}=n_{2}=2$) is a double (2,2)-form which is 
symmetric (the pairs can be interchanged); the Ricci tensor  
is a double symmetric (1,1)-form and, in general, any completely symmetric 
rank-$r$ tensor is an $r$-fold (1,1,\dots ,1)-form. 
A 3-tensor $A_{abc}$ with the 
property $A_{abc}=-A_{cba}$ is a double 
(2,1)-form and the corresponding $\tilde A$ is clearly given by 
$\tilde{A}_{abc}=\tilde{A}_{[ab]c}\equiv A_{acb}$. 
The standard index version of familiar tensors such as Riemann tensors
$R_{abcd}=R_{[ab][cd]}$, Weyl tensors $C_{abcd}=C_{[ab][cd]}$, torsion
tensors $T_{abc}=T_{[ab]c}$, or Lanczos tensors $H_{abc}=H_{[ab]c}$,
already have the indices in the appropriate permuted version, so 
that they coincide with their tilded versions. In these cases
we shall dispense with the \ $\tilde{}$ \ label.
%% 
 % and shall always
 % consider the tensors which we are using to have undergone the required
 % permutation of indices and to be in the standarised $r$-fold form
 % described above.
 %%

The scalar product (\ref{scalar-product}) can be immediately 
generalized to arbitrary tensors of the same type (i.e., with 
the same form-structure number $r$ and block ranks) by first defining
$$
(T,S) \equiv \tilde{T}_{a_{1}\dots a_{m}}\tilde{S}^{a_{1}\dots a_{m}}
$$
and then (same comments as before on compactness)
\be
<T , S > \, \equiv \int_{V_{n}} (T,S) \, \eta \, .\label{scalar-p}
\ee
Observe that this scalar product is adapted to the structure as 
$r$-fold forms of the tensor fields $T$ and $S$, and therefore 
it is different from the product defined in (\ref{scalarL}), so that
we have, in general,  $<T , S >\, =\, \{\tilde T , \tilde S \}\, \neq \{T,S\}$. 
This will be relevant in what follows. Note, moreover, that $\{\, ,\, 
\}$ is defined for general tensors $T,S$, not only for those with the 
same form-structure number $r$ and block ranks. 
As before, the product (\ref{scalar-p}) is bi-linear, symmetric and 
non-degenerate. For proper Riemannian manifolds, it is also positive 
definite.

Now consider any tensor field $T$ and let $r$ be its form-structure number. 
For each one of the $r$ antisymmetric blocks, one can 
follow a procedure similar to that recalled above for $p$-forms. 
To fix ideas, 
let us select the $i$-th block, with block rank $n_{i}$, to 
treat $T$ as a tensor-valued 
$n_{i}$-form, so that 
 $\tilde{T}^{a_1\ldots a_h} {}_{b_1\dots b_{n_i}}{}^{a_{h+1}\ldots a_k}=
\tilde{T}^{a_1\ldots a_h}{}_{[b_1\dots b_{n_i}]}{}^{a_{h+1}\ldots a_k}$ with  
a set of $n_i\leq n$ {\it completely antisymmetrical} indices
${b_1\dots b_{n_i}}$ plus a number of extra indices 
denoted by ${a_1\ldots a_h}$ and ${a_{h+1}\ldots a_k}$. 
These objects are usually referred to as {\it differential tensor forms} 
or {\it tensor-valued differential forms}, e.g. \cite{Car,Bel,MTW,Coq}. 
Then, the previous definition (\ref{*}) can be trivially extended so that 
the dual of $T$ with respect to the $i$-th block $b_1\dots b_{n_i}$
is denoted as $*_{(i)} T$ and defined by \cite{S}
\be
(*_{(i)} T)^{a_1\ldots a_h}{}_{b_{n_i+1}\dots b_n}{}^{a_{h+1}\ldots a_k} 
%\tilde{T}^{a_1\ldots a_h}{}_{\stackrel{*}{b_{n_i+1}\dots b_n}}{}^{a_{h+1}\ldots a_k}=
%T^{a_1\ldots a_h}{}_{\stackrel{*}{[b_{n_i+1}\dots b_n]}}{}^{a_{h+1}\ldots a_k}
\equiv
\frac{1}{{n_i}!}\eta_{b_1\dots b_n}\tilde{T}^{a_1\ldots a_h}{}^{b_1\dots 
b_{n_i}}{}^{a_{h+1}\ldots a_k} \, .
\label{*2}
\ee
%where the * is placed over the indices onto which the dual operation acts.
As before, $*_{(i)}*_{(i)}=\epsilon (-1)^{n_{i}(n-n_{i})}$ when acting on $\tilde T$.

Next we extend the definitions of $d, \delta, \D$ from the previous
section to differential $d_{(i)}$, co-differential
$\delta_{(i)}$, and de Rham operator $\D_{(i)}$ acting on the 
$i$-th block of the $r$-fold form $T$ as:\footnote{When applying the 
operators $*_{(i)},d_{(i)},\delta_{(i)},\D_{(i)}$ we are defining the 
resulting tensor as having the $i$-th block in the ordered position, 
so that  $\widetilde{*_{(i)}T}=*_{(i)}T$, $\widetilde{d_{(i)}T}=d_{(i)}T$, 
etc.}
\bea
(d_{(i)} T)^{a_1\ldots a_h}{}_{b_{1}\dots b_{n_{i}+1}}{}^{a_{h+1}\ldots a_k} &\equiv&
(-1)^{n_{i}}(n_{i}+1)
\tilde{T}^{a_1\ldots a_h}{}_{[b_{1}\dots b_{n_{i}}}{}^{a_{h+1}\ldots a_k}{}_{;b_{n_{i}+1}]}\, ,\label{newd} \\
\delta_{(i)} = (-1)^{n_{i}} *_{(i)}^{-1}d_{(i)} *_{(i)} &\equiv& 
\epsilon\, (-1)^{(n-n_{i})(n_{i}-1)+1} *_{(i)} d_{(i)} *_{(i)} \label{newdelta0}\\
\D_{(i)}& \equiv& d_{(i)}\delta_{(i)}+\delta_{(i)}d_{(i)}\, .\hspace{2cm}\label{newDelta}
\eea
(An alternative equivalent definition for $\delta_{(i)}$ can be given via $d_{(i')}$ and a contraction (trace) operator across the $i$-th and $i'$-th blocks \cite{Kul}; see comments following the definition of the trace operator in Section 5.)

We emphasise that the covariant derivative acts on {\it all} indices,
i.e., the extra tensor indices as well as the explicit form indices. 
Although the definition for $d$ was extended to a tensor valued form
in \cite{Car} and is well known, as far as we are aware {\it all three
operators}, in the form given above, were first introduced in
\cite{Bel}, and are less familiar.  In this more general context,
sometimes the extended $d$ is called `the absolute exterior
differential' or `the covariant exterior differential'; some
references retain $d$ and $\delta$ \cite{MTW}, some texts replace 
$d\ (\delta ,\Delta)$ with $\tilde d \ (\tilde {\delta},\tilde{\Delta})$
\cite{Bel}, or with $d^{\nabla} \ (\delta^{\nabla},\Delta^{\nabla})$
\cite{Bes,Coq,Ur}, or with $D$ \cite{Car,Kul,Th,BT,GS,Na,BCJR,BCJF}, or
with $\nabla$ \cite{Fr}.  We shall use the notation of (\ref{newd},
\ref{newdelta0}, \ref{newDelta}) since we wish to be able to combine
different de Rham operators.

The operator $d_{(i)}$ ($\delta_{(i)}$) produces another tensor 
with one more (one less) index in general, and this new tensor has
in general the same form-structure number $r$. 
Of course, as for single $p$-forms, there are some special situations: 
(i) if $n_{i}=n$ then $d_{(i)}T=0$; (ii) for any tensor $T$ with form 
structure number $r$, $\nabla 
T$ has form-structure number $r+1$ and can be considered as a 
definition of `$d_{(r+1)}T$' ; (iii) if $n_{i}=1$ then
$\delta_{(i)}T$ has $r-1$ as form-structure number. In this case, in 
order to compute $\D_{(i)}$ one 
has to allow the operator $d_{(i)}$ in the combination 
$d_{(i)}\delta_{(i)}$ to act on the missing block---as if 
$\delta_{(i)}T$ were an $r$-fold 
$(n_{1},\dots ,n_{i-1},0,n_{i+1}\ldots n_r)$-form--- in the same way 
as in (ii). 
%% 
 % In this way, there is no difficulty when we are considering $\D_{(i)}$, since the 
 % operator $d_{(i)}$ can act on the $i$-th block, as required.
 %%

Taking this into account, we need definitions for 
$d_{(\ell)}T,\delta_{(\ell)}T$ 
when $\ell\neq 1,\dots ,r$, which are simply
\be
d_{(\ell)}T = \widehat{\nabla \tilde{T}} ,  \hspace{1cm}
\delta_{(\ell)} T =0, \hspace{1cm} \ell \notin \{1,\dots ,r\}
\label{extra}
\ee
where in the first case the $\ \widehat{}\ $ means that the extra index 
provided by the covariant derivative must be placed in the appropriate 
place within $\{1,\dots ,r+1\}$.

Bearing all this in mind, we can write down explicit index formulas for 
$\delta_{(i)}$, $\D_{(\ell)}$ and $\D_{(i)}$ for the cases $i\in 
\{1,\dots ,r\}\not\ni \ell$ :
\bea
&&(\delta_{(i)}  T)^{a_1\ldots a_h}{}_{b_{1}\dots b_{n_{i}-1}}{}^{a_{h+1}\ldots a_k}= 
-\, \tilde{T}^{a_1\ldots a_h}{}{}_{cb_{1}\dots b_{n_{i}-1}}{}^{a_{h+1}\ldots a_k}{}^{;c}\, ,
\hspace{1cm} \label{newdelta} \\
&& (\D_{(\ell)} T)_{a_{1}\dots a_{m}} = - \nabla^c\nabla_{c}
\tilde{T}_{a_{1}\dots a_{m}}
\label{trivialLap} \\
&&(\D_{(i)} T)^{a_{1}\dots a_{h}}{}_{b_{1}\dots b_{n_{i}}}{}^{a_{h+1}\ldots a_k}=-
\nabla^c\nabla_{c} \tilde{T}^{a_{1}\dots a_{h}}{}_{b_{1}\dots b_{n_{i}}}{}^{a_{h+1}\ldots a_k}\nonumber \\ && \quad+
n_{i}R_{c[b_{1}}\tilde{T}^{a_{1}\dots a_{h}c}{}_{b_{2}\dots b_{n_{i}}]}{}^{a_{h+1}\ldots a_k}
 -\, \frac{n_{i}(n_{i}-1)}{2}R_{cd[b_{1}b_{2}}
\tilde{T}^{a_{1}\dots a_{h}cd}{}_{b_{3}\dots b_{n_{i}}]}{}^{a_{h+1}\ldots a_k}\nonumber \\
&&\qquad -
n_{i}\sum_{s=1}^k R_{c}{}^{a_{i}}{}_{d[b_{1}}\, 
\tilde{T}^{a_{s}\dots a_{s-1}\,c\, a_{s+1}\dots a_{h}d}{}_{b_{2}\dots 
b_{n_{i}}]}{}^{a_{h+1}\ldots a_k}\, .\label{newDelta2}
\eea

A trivial calculation using the Gauss theorem leads to
$$
<d_{(i)}T , U > \, = \, <T , \delta_{(i)} U >
$$
where $U$ and $T$ have the same form-structure number $r$ and 
$n_{i}(T)+1=n_{i}(U)$. Thus, the operators $d_{(i)}$ and $\delta_{(i)}$ 
are adjoint to each other with respect to $<\, ,\, >$ introduced in 
(\ref{scalar-p}). 

Observe the similarities and the differences of $\D_{(i)}$ of 
(\ref{newDelta2}) 
with $\D_{L}$ as defined 
in (\ref{Li}): 
\begin{itemize}
\item For each $i=1,\dots ,r$, $\D_{(i)}T$ respects the skew-symmetry on the 
$i$-th antisymmetric block of $T$, and the symmetries and trace properties 
on the extra indices not in that block; this implies, in particular, that 
$\D_{(i)}T$ has the same form-structure number and block ranks as $T$. 
\item On the other hand, any {\em mixed trace}, or {\em 
mixed index symmetry}, involving indices from {\it both} the explicit 
$i$-th antisymmetric 
block and the rest of the indices is {\em not} preserved in general.
\item For each $i=1,\dots ,r$, $\Delta_{(i)}$ is self-adjoint with respect to 
the scalar product (\ref{scalar-p}), so that $<\D_{(i)}T,S>\, =\, 
<T,\D_{(i)}S>$ for all $T,S$ with the same form-structure number and block ranks.
{\em Moreover}, one can prove the identities
\be
<T ,\D_{(i)}T >\, =\, <d_{(i)}T ,d_{(i)}T > + <\delta_{(i)}T ,\delta_{(i)}T> 
\hspace{3mm} \forall i\in\{1,\dots ,r\} \label{i-identity}
\ee 
in the same way as (\ref{dotp}) was obtained.  From here, in particular, one deduces that $\D_{(i)}$ are positive 
operators in proper Riemannian manifolds: $<T, \D_{(i)} T> \,\, \geq 0$.
\item When $r=1$, that is to say, when acting on $p$-forms, $\D_{(1)}$ 
coincides with the de Rham operator: $\D_{(1)}\S=\D\S=\D_{L}\S$ for all 
$\S\in \Lambda^p$. In particular, when acting on rank-1 tensors, $\D_{(1)}$ 
commutes with the 
covariant derivative in manifolds with a parallel Ricci tensor ($\nabla_aR_{ab}=0$). There is however no simple commutation relation between $\D_{(i)}$ and $\delta_{(i)}$; see (\ref{comm2}) below.
\end{itemize}

As with $p$-forms, we can put forward the following definition:
\begin{defi}[$i$-harmonic and fully harmonic tensors]
A tensor field $T$ with 
form-structure number $r$ is said to be {\em i-harmonic}, for $i\in\{1,\dots 
,r\}$, if and only if $\D_{(i)}T=0$. Such a tensor will be called {\em 
fully harmonic} if it is $i$-harmonic {\em for all} 
$i=1,\dots ,r$. 
\label{i-harmonic}
\end{defi}
Note that the harmonic property in the sense of Lichnerowicz
(i.e., $\D_{L}T=0$) is different from these new harmonic properties. 
This may have important consequences for the theory of harmonic 
tensors in Riemannian or pseudo-Riemannian manifolds. See Definition 
\ref{harmonic} and the  results which follow it.

It is important to note, unlike the situation for $p$-forms, that  
$d_{(i)}^2\ne 0$ and $\delta_{(i)}^2 \ne 0$ in curved spaces. 
Using the Ricci identity, one obtains 
\bea
& & (d_{(i)}^2{ T})^{a_{1}\dots a_{h}}{}_{b_{1}\dots b_{n_{i}+2}}{}^{a_{h+1}\ldots a_k} \nonumber \\ & & \qquad =
\frac{1}{2}(n_{i}+1)(n_{i}+2)\sum_{s=1}^k 
R^{a_{s}}{}_{c[b_{n_{i}+1}b_{n_{i}+2}}\,
{\tilde{T}}^{a_{1}\dots a_{s-1}ca_{s+1}\dots a_{h}}{}_{b_{1}\dots b_{n_{i}}]}{}^{a_{h+1}\ldots a_k}\, ,
\label{dr2}\\
& & (\delta_{(i)}^2{T})^{a_{1}\dots a_{h}}{}_{b_{1}\dots b_{n_{i}-2}}{}^{a_{h+1}\ldots a_k}=
-\frac{1}{2}\sum_{s=1}^k 
R^{a_{s}}{}_{c}{}^{de}\,
{\tilde{T}}^{a_{1}\dots a_{s-1}ca_{s+1}\dots a_{h}}{}_{deb_{1}\dots b_{n_{i}-2}}{}^{a_{h+1}\ldots a_k}\, .
\label{deltar2}
\eea
Note that for higher derivatives all derivatives of the Riemann tensor disappear
because of the Bianchi identities; this means that even derivatives will give terms involving products of the Riemann tensor and $\tilde T$, while odd derivatives will give terms involving products of the Riemann tensor and a derivative of $\tilde T$.

have a similar product structure with the
Riemann tensor, since all derivatives of the Riemann tensor disappear
because of the Bianchi identities.

Observe that, in flat spaces, these operators are nilpotent. (Compare 
with \cite{DH}, where similar ideas are developed in flat Euclidean 
spaces.) Note also that $d_{(i)}^2=0$ whenever $n_{i}\in \{n-1,n\}$, and 
that $\delta_{(i)}^2=0$ if $n_{i}=1$. More generally, $d_{(i)}^m =0$ for 
$m>n-n_i$, and $\delta^m_{(i)}=0$ for $m>n_i$, so these curved space operators 
are also nilpotent, but depending on the dimension of the space and the block 
rank of the tensor-valued form being acted on.

The 
commutation properties with $\D_{(i)}$ follow directly from 
(\ref{dr2}-\ref{deltar2}) using
\be
d_{(i)}\D_{(i)}-\D_{(i)}d_{(i)}=d_{(i)}^2\delta_{(i)}-\delta_{(i)}d_{(i)}^2 \, , 
\hspace{1cm}
\delta_{(i)}\D_{(i)}-\D_{(i)}\delta_{(i)}=\delta_{(i)}^2d_{(i)}-d_{(i)}\delta_{(i)}^2 
\, . \label{comm2}
\ee
(Observe that, in flat spaces, these operators commute.) Finally, 
$\D_{(i)}$ commute with each of the Hodge dual operators $*_{(j)}$
defined in (\ref{*2})
$$
*_{(j)}\,\, \D_{(i)} = \D_{(i)}\, *_{(j)} \hspace{1cm} 
\forall i,j\in \{1,\dots ,r\}\, .
$$
%where the $*_{(j)}$ dualizes the $j$-th block.

The usefulness of these operators (\ref{newd}), (\ref{newdelta0}),
(\ref{newDelta}) is that we can now translate the calculation
performed in (\ref{calcul}) to arbitrary tensors.  Although the
operators $\D_{(i)}$ are different from the de Rham operator $\D$ it is
still possible  to apply the Cauchy-Kovalewski theorem
\cite{CH}, and we can generalise Result \ref{r2} to
\begin{result}
    \label{1}
Let $T$ be any tensor field, and let $r$ be its form-structure number 
and $(n_{1},\dots ,n_{r})$ its block ranks. There always exist local 
superpotentials $\stackrel{o}{\mbox{$T$}}$ with the same form-structure number 
and block ranks as $T$ such that $\D_{(i)} \stackrel{o}{\mbox{$T$}}=T$; 
furthermore, for each $i=1,\dots ,r$ there always exist a pair of {\em local}
potentials $(Y_{(i)},Z_{(i)})$, such that
$$
T=\delta_{(i)}(d_{(i)}\stackrel{o}{\mbox{$T$}})+
d_{(i)}(\delta_{(i)}\stackrel{o}{\mbox{$T$}})\equiv \delta_{(i)}Y_{(i)} 
+d_{(i)}Z_{(i)}
$$
where  we have defined the $r$-fold form potentials 
$Y_{(i)}=d_{(i)}\stackrel{o}{\mbox{$T$}}$ and 
$Z_{(i)}=\delta_{(i)}\stackrel{o}{\mbox{$T$}}$.   
\end{result}

In index notation:
$$
(\D_{(i)}{ \stackrel{o}{\mbox{$T$}}}){}^{a_1\ldots a_h}{}_{b_{1}\dots b_{n_{i}}}{}^{a_{h+1}\ldots a_k}=
\tilde{T}^{a_1\ldots a_h}{}_{b_{1}\dots b_{n_{i}}}{}^{a_{h+1}\ldots a_k}
$$
and
\be
\tilde{T}^{a_1\ldots a_h}{}_{b_{1}\dots b_{n_{i}}}{}^{a_{h+1}\ldots a_k}
=- Y^{a_1\ldots a_h}{}_{cb_{1}\dots b_{n_{i}}}{}^{a_{h+1}\ldots a_k}{}^{;c}+
(-1)^{n_{i}-1}n_{i}\,
Z^{a_1\ldots a_h}{}_{[b_{1}\dots b_{n_{i}-1}}{}^{a_{h+1}\ldots a_k}{}_{;b_{n_{i}}]}
\label{twopot}
\ee
where $Y^{a_1\ldots a_h}{}_{b_{1}\dots b_{n_{i}+1}}{}^{a_{h+1}\ldots a_k}$ is a 
rank-$(m+1)$ tensor of 
type $r$-fold $(n_{1},\dots ,n_{i}+1,\dots ,n_{r})$-form and
$Z^{a_1\ldots a_h}{}_{b_{1}\dots b_{n_{i}-1}}{}^{a_{h+1}\ldots a_k}$ is a 
rank-$(m-1)$ tensor of 
type  $r$-fold $(n_{1},\dots ,n_{i}-1,\dots ,n_{r})$-form, 
given respectively by 
\bea
Y^{a_1\ldots a_h}{}_{b_{1}\dots b_{n_{i}+1}}{}^{a_{h+1}\ldots a_k}&\equiv& (-1)^{n_{i}}
(n_{i}+1){\stackrel{o}{\mbox{$T$}}}{}^{a_1\ldots a_h}{}_{[b_{1}\dots 
b_{n_{i}}}{}^{a_{h+1}\ldots a_k}{}_{;b_{n_{i}+1}]} \, , \nonumber\\  Z^{a_1\ldots a_h}{}_{b_{1}\dots b_{n_{i}-1}}&\equiv& -\, 
{\stackrel{o}{ \mbox{$T$}}}{}^{a_1\ldots a_h}{}_{cb_{1}\dots b_{n_{i}-1}}{}^{a_{h+1}\ldots a_k}{}^{;c}
\eea 

Trivially, these potentials in the above result are highly non-unique. 
For $r$-fold forms with $r>1$ we can repeat the above constructions for 
each $i\in \{1,\dots ,r\}$. 
Therefore the previous definition of potentials
can be applied to each of the 
antisymmetric blocks, and thus 
there are many possible different pairs of potentials. 

A second disadvantage is that each of the operators $\D_{(i)}, \, i=1,2,
\ldots r$, in general, does not have all of the useful properties
which we would like.  Each operator is adapted to respect the index
properties of its corresponding block, but not for the other blocks of
indices.  So, for instance, any index symmetry or trace property
across the different blocks of a tensor will not necessarily be directly
reflected in its superpotential.

\section{A weighted de Rham operator and associated potentials. 
Harmonic tensors.}

Of course, one can mix the operators $\D_{(i)}, \, i=1,2, \ldots r$,
weighting them, and obtain new Laplace-type operators with similar or
better properties.  A particularly good one is,
\begin{theorem}\label{barD}
The operator $\bar\D$ given by
\be
\bar{\D}\equiv \frac{1}{r} (\D_{(1)}+\D_{(2)}+\dots +\D_{(r)}) \label{av-Delta}
= \frac{1}{r}\sum_{i=1}^r \D_{(i)} \ee
is linear, self-adjoint with respect to the scalar product 
(\ref{scalar-p}), respects all index symmetry 
properties, and commutes with all trace operations when acting on 
arbitrary tensor fields (here $r$ is the form-structure number of the 
tensor field). It is related to the Lichnerowicz operator by 
\be
\, \bar{\D}=\frac{1}{r}\Delta_{L}-\frac{r-1}{r}\nabla^c\nabla_{c} \ . 
\label{relation}
\ee
\end{theorem}
\begin{proof}
The linearity and self-adjointness follows from these properties for 
each of the $\D_{(i)}$. The formula (\ref{relation}) can be easily 
established from (\ref{Li}) and (\ref{newDelta2}) (see also formulas 
(\ref{barDelta}) and (\ref{barDelta11}-\ref{barDelta24}) below). 
The commutativity with index permutations and 
traces then follows immediately from the properties of $\Delta_L$ and 
$\nabla^c\nabla_{c}$. \square
\end{proof}

So we find that although $\D_L$ and $\bar \D$  do not coincide, they are closely related, and indeed $\bar \D$ has most of the useful properties of  $\D_L$; but crucially  $\bar \D$ in addition has direct links with $d_{(i)}$ and $\delta_{(i)}$.
Hence we believe that $\bar\D$ is a  more powerful alternative than the classical Lichnerowicz
operator $\Delta_L$  since it is so well adapted to dealing
with the $r$-fold form structure of the tensors. 
%% 
 % note that, like $\Delta_L$, $\bar \Delta$ is
 % self-adjoint, but it does not have the commutativity properties shown
 % for $\Delta_L$ in \cite{Lic}.  
 %%
Of course, for single $p$-forms we have
$$
\bar\D \S=\D \S=\D_{L} \S , \  \  \forall \S\in \Lambda^p \, .
$$

An extremely important consequence of this new operator $\bar\D$
defined in (\ref{av-Delta}) is that for any given tensor field
$T$, there exists an {\it associated}
superpotential $\stackrel{o}{\mbox{$T$}}$, by which we
mean a superpotential, not just with the same form-structure number 
and block ranks, but also {\it with
the same index symmetries and trace properties} as $T$.  
Moreover, from this associated superpotential, a set of $2r$
potentials can be obtained in a natural and straightforward manner 
from the definition of $\bar\D$
$$
\bar\D \stackrel{o}{\mbox{$T$}} \ \equiv \frac{1}{r}\sum_{i=1}^r
[\delta_{(i)}(d_{(i)}\stackrel{o}{\mbox{$T$}})+
d_{(i)}(\delta_{(i)}\stackrel{o}{\mbox{$T$}})] \, .
$$
Explicitly we have thus proven, 
\begin{theorem}
    \label{3}
Given any  tensor field $T$, there always
exists an {\em associated} superpotential $\stackrel{o}{\mbox{$T$}}$ 
such that $\bar\D \stackrel{o}{\mbox{$T$}} = T$; 
furthermore, if $r$ is the form-structure number of $T$,
there always exists a set of $2r$
{\em local} potentials $(Y_{(i)},Z_{(i)})$, $i=1,2,\ldots r$, such
that
$$
T=\frac{1}{r}\sum_{i=1}^r (\delta_{(i)}Y_{(i)} +d_{(i)}Z_{(i)})
$$
where $Y_{(i)}=d_{(i)}\stackrel{o}{\mbox{$T$}}$ and
$Z_{(i)}=\delta_{(i)}\stackrel{o}{\mbox{$T$}}$ are the potentials.
\end{theorem}
In our earlier derivation of the potential for the Weyl tensor
\cite{ES} the operator which we constructed for the superpotential is
the special case of $\bar\D$ for a double $(2,2)$-form.  To the best
of our knowledge, this seems to be the only time such a {\it weighted
de Rham operator} $\bar\D$ has appeared in the literature.  There are
a variety of other generalised Laplacians \cite{Bes,Ur,Na}, 
but for tensor-valued forms in pseudo-Riemannian spaces
with Levi-Civita connections the only Laplacians mentioned are usually
$-\nabla^c\nabla_{c}$, Lichnerowicz's Laplacian $\D_L$ and the Laplacians
which we have labeled $\D_{(i)}$.\footnote{To avoid any
misunderstanding we point out that what we have called the `usual
Laplacian' $-\nabla^c\nabla_{c}$ (and with the symbol $\D_{(\ell)}$ in (\ref {trivialLap})) is sometimes given by the symbol $\bar\D$,
\cite{Lic,Ur} and called the `rough Laplacian' \cite{Ur} or 
the `naive Laplacian' \cite{Coq}.}  In
particular for double $(2,2)$-forms such as the Riemann tensor, the
Laplacian used in \cite{Bel,MTW,BF,BCJF,BCJR}, is
the one we label $\D_{(2)}$, and as we argue in Section 7, this is not
really a suitable operator in that context.  It is only in \cite{Coq}
that any attempt is made to compare different Laplacians, where it is
pointed out that, for the special case of a symmetric $2$-tensor, 
the Lichnerowicz Laplacian can be written as 
$\D_L = \D_{(1)}+ \D_{(2)} +\nabla^c\nabla_{c}$ (which can
be considered as a special case of formula (\ref{relation})); but the
significance of $(\D_{(1)}+\D_{(2)})$ as a Laplacian was not recognised in
\cite{Coq}.

We start with a simple logical definition which will allow us to show 
the potentialities of the new Laplacian $\bar\D$ with a very simple 
preliminary application.
\begin{defi}
A tensor field $T$ will be called {\em harmonic} if and only if \label{harmonic}
$$
\bar\D T=0 \, .
$$\end{defi}
Obviously, any fully harmonic tensor in the sense of Definition 
\ref{i-harmonic} is trivially harmonic. The converse, however, does 
not hold in general. Nevertheless, from formula (\ref{i-identity}) we deduce 
the general identity
\be
<T ,\bar\D T >\, =\frac{1}{r} \sum_{i=1}^r <T , \D_{(i)} T > \, =\frac{1}{r}
\sum_{i=1}^r\left(<d_{(i)}T ,d_{(i)}T > + <\delta_{(i)}T ,\delta_{(i)}T> 
\right) \label{bar-identity}
\ee
for arbitrary tensor fields. Therefore, it is straightforward to obtain the 
following converse in proper Riemannian manifolds.
\begin{theorem}
\label{harmonic=full}
Let $(V_{n},g)$ be a compact without boundary proper Riemannian 
manifold. Then, a tensor $T$ is harmonic if and only if it is fully 
harmonic, and if and only if 
$$
d_{(i)}T=0, \hspace{1cm} \delta_{(i)}T=0,\hspace{1cm} \forall 
i\in\{1,\dots,r\}\, .
$$
\end{theorem}
{\bf Remarks}
\begin{itemize}
\item
 The above definition is different from the Lichnerowicz 
harmonic property which makes use of the condition $\D_{L}T=0$; and we emphasise
 that, in general,  there is no analogue of
this theorem by using $\D_{L}$.

\item This definition differs from the harmonic definition in \cite{Kul} for $2$-fold forms in proper Riemannian space;  in fact the definition there coincides with our definition of {\it $1$-harmonic} $2$-fold forms, which is equivalent to $\Delta_1 \omega=0$ where $\omega $ is a $(p,q)$-form.

\item See Section \ref{sec:Riemann} for details concerning the harmonic 
property of the Riemann tensor. 

\item The tensor decomposition into
superpotentials and potentials for vectors and $2$-tensors developed
in \cite{Gem} exploited the Lichnerowicz operator $\Delta_L$ in order
to find the superpotential $U$ for any tensor $V$ via $\Delta_L U = V$;
however, such a superpotential $U$ does not, in general, yield
potentials in a straightforward manner as there is no direct natural link to
any first order operators.
\end{itemize}

 For general calculations, we would need the
commutation properties of the operators $d_{(i)}$ and $\delta_{(j)}$. 
Actually the general case requires a detailed and cumbersome study in
non-flat backgrounds; but we are only interested in the case of {\it
double forms} ($r=2$) in this paper, and so we devote the next section to this
case in full detail.

\section{Double $(q,p)$-forms}
In this section we deal with tensors with form-structure number $r=2$ 
and block ranks $n_{1}=q$ and $n_{2}=p$, called double $(q,p)$-forms. These 
are tensors with $p+q$ indices which fall into two categories: 
a block of $q$ antisymmetric indices, and another block of $p$ 
antisymmetric indices. By rearranging if necessary,
the tilded version of the tensor can then always be written as
$$
\tilde{T}^{a_{1}\dots a_{q}}{}_{b_{1}\dots b_{p}}=
\tilde{T}^{[a_{1}\dots a_{q}]}{}_{[b_{1}\dots b_{p}]}.
$$
In addition to the Lichnerowicz operator, there are three other
Laplace-like operators $\D_{(1)}$, $\D_{(2)}$ and $\bar\D$ as defined 
in the previous sections acting on these objects. Using (\ref{newDelta}), 
(\ref{newDelta2}), and (\ref{av-Delta}) their explicit formulas are
\bea
(\D_{(1)}T)^{a_{1}\dots a_{q}}{}_{b_{1}\dots b_{p}}= 
-\nabla^c\nabla_{c} \tilde{T}^{a_{1}\dots a_{q}}{}_{b_{1}\dots b_{p}}+
qR^{c[a_{1}}\tilde{T}_{c}{}^{a_{2}\dots a_{q}]}{}_{b_{1}\dots b_{p}}\\ \nonumber
-\frac{q(q-1)}{2}
R^{cd[a_{1}a_{2}}\tilde{T}_{cd}{}^{a_{3}\dots a_{q}]}{}_{b_{1}\dots b_{p}}+
qp\, R^{[a_{1}}{}_{cd[b_{1}} \tilde{T}^{a_{2}\dots a_{q}]cd}{}_{b_{2}\dots b_{p}]} \, ;\label{barDelta1T}
\eea
\bea
(\D_{(2)}T)^{a_{1}\dots a_{q}}{}_{b_{1}\dots b_{p}}=
- \nabla^c\nabla_{c} \tilde{T}^{a_{1}\dots a_{q}}{}_{b_{1}\dots b_{p}}+
pR_{c[b_{1}}\tilde{T}^{a_{1}\dots a_{q}c}{}_{b_{2}\dots b_{p}]}\\ \nonumber
-\frac{p(p-1)}{2} R_{cd[b_{1}b_{2}} 
\tilde{T}^{a_{1}\dots a_{q}cd}{}_{b_{3}\dots b_{p}]}
+pqR^{[a_{1}}{}_{cd[b_{1}} \tilde{T}^{a_{2}\dots a_{q}]cd}{}_{b_{2}\dots b_{p}]} \, ;\label{barDelta2T}
\eea
\bea
(\bar\D T)^{a_{1}\dots a_{q}}{}_{b_{1}\dots b_{p}}=
- \nabla^c\nabla_{c} \tilde{T}^{a_{1}\dots a_{q}}{}_{b_{1}\dots b_{p}}
+\frac{q}{2}R^{c[a_{1}}\tilde{T}_{c}{}^{a_{2}\dots a_{q}]}{}_{b_{1}\dots b_{p}}
+\frac{p}{2}R_{c[b_{1}}\tilde{T}^{a_{1}\dots a_{q}c}{}_{b_{2}\dots b_{p}]}
 \nonumber \\
-
\frac{q(q-1)}{4}
R^{cd[a_{1}a_{2}}\tilde{T}_{cd}{}^{a_{3}\dots a_{q}]}{}_{b_{1}\dots b_{p}}-
\frac{p(p-1)}{4} R_{cd[b_{1}b_{2}} \tilde{T}^{a_{1}\dots a_{q}cd}{}_{b_{3}\dots b_{p}]}
 \nonumber \\+qpR^{[a_{1}}{}_{cd[b_{1}} \tilde{T}^{a_{2}\dots a_{q}]cd}{}_{b_{2}\dots b_{p}]} 
\label{barDelta}
\eea
The same comments as before for the cases with $p$ or $q$ equal to 
$1$, $n-1$ or $n$ are in order here.
We also note from formula (\ref{relation}) in Theorem \ref{barD} that  
$\bar \D$ is related to the Lichnerowicz operator in this case by 
$$
\bar\D=\frac{1}{2}(\D_{L}-\nabla_c\nabla^c) .
$$

Specialising Theorem \ref{3} we obtain,
\begin{coro}
Given any tensor field $T$ with the structure of a double $(q,p)$-form
there always exists an {\em associated} local superpotential 
${\stackrel{o}{\mbox{$T$}}}$ such that
$\bar\D \stackrel{o}{\mbox{$T$}} =T$; furthermore there always exist 
local potentials $Y_{(1)}, Y_{(2)}, Z_{(1)}, Z_{(2)}$ such that
\be
T=\frac{1}{2}\left(\delta_{(1)}Y_{(1)}+
\delta_{(2)}Y_{(2)}+d_{(1)}Z_{(1)}+d_{(2)}Z_{(2)}\right),\label{4pot}
\ee
where the double $(q+1,p)$-form $Y_{(1)}=d_{(1)}\stackrel{o}{\mbox{$T$}}$, 
the double $(q,p+1)$-form $Y_{(2)}=d_{(2)}\stackrel{o}{\mbox{$T$}}$, 
the double $(q-1,p)$-form $Z_{(1)}=\delta_{(1)}\stackrel{o}{\mbox{$T$}}$ 
and the double $(q,p-1)$-form $Z_{(2)}=\delta_{(2)}\stackrel{o}{\mbox{$T$}}$.
\end{coro}

In the sequel we will need to deal with contractions as well as permutations of indices for different
expressions, so we make two formal definitions which will be useful.
\begin{defi} [Trace of a double form]
The {\em trace} of a double
$(q,p)$-form $T$ is the double $(q-1,p-1)$-form $\hbox{\em tr}(T)$ given
by
\be 
(\widetilde{\hbox{\em tr}(T)})^{a_2\ldots a_q}{}_{b_2 \ldots b_p} \equiv 
\tilde T^{ca_2\ldots a_q}{}_{cb_2 \ldots b_p} \, .\label{trdefn}
\ee
\end{defi}
Once more, we must remark that if $q=1$ (or $p=1$), then the 
first (second) block disappears after taking the trace tr, so that the 
resulting tensor has a form-structure number less than 2. In these 
situations, and as happened in other situations explained before, sometimes 
it is necessary to consider the resulting tensor tr$(T)$ as an 
equivalent double $(0,p-1)$-form (or double 
$(q-1,0)$-form). Bearing 
this in mind, some useful results follow immediately from the respective 
definitions 
(\ref{newd}) and (\ref{newdelta}),
\be
\hbox{tr}(d_{(1)} T) =   -d_{(1)} \hbox{tr} (T)    -\delta_{(2)} T ,\hspace{1cm}
\hbox{tr}(d_{(2)} T) =   -d_{(2)} \hbox{tr} (T)    -\delta_{(1)} T , \label{trd}
\ee
\be 
\hbox{tr}(\delta_{(1)} T) =  -\delta_{(1)} \hbox{tr} (T) \,\,\, (q\ge 2); 
\hspace{1cm}
\hbox{tr}(\delta_{(2)} T) =  -\delta_{(2)} \hbox{tr} (T) \,\,\, (p\ge 2). 
\label{trdelta}
\ee
In \cite{Kul} the second of equations (\ref{trd}) was used as the definition of $\delta_{(1)}$ for a double $(p,q)$-form.
\begin{defi} [Transpose of a double form]
For a double $(q,p)$-form $T$ we  define the
{\em (generalised) transpose} ${}^t T$ of $T$ as the double 
$(p,q)$-form given by interchange of the blocks:
$$
({\widetilde{{}^t T}})^{a_{1}\dots a_{p}}{}_{b_{1}\dots b_{q}} \equiv
\tilde{T}_{b_{1}\dots b_{q}}{}^{a_{1}\dots a_{p}}. \label{trandef}
$$
\end{defi}
Obviously ${}^{tt} T = T$, and the following useful results also follow from the definitions
\be
d_{(2)}({}^tT)=  {}^t (d_{(1)}T) ,\hspace{1cm} \delta_{(2)}({}^t T) =\  {}^t (\delta_{(1)}T)
\label{trans-ddel}
\ee
In addition, we note that the operations of transpose and trace commute, $${}^t(\hbox{tr}(T))= \hbox{tr}({}^tT) \ .$$

Of course, we could also define traces  and transposes
of the more general $r$-fold forms in the previous section, 
by taking pairs of blocks at a time.

\

For future reference we give explicitly the five cases of $(1,1)$-,
$(2,1)$-, $(2,2)$-, $(2,3)$- and $(2,4)$-forms respectively for
(\ref{barDelta})
\bea
-(\bar\D T)^{a_{1} }{}_{b_{1}}=
\nabla^c\nabla_{c} T^{a_{1} }{}_{b_{1}}
-{1\over 2}R^{ca_{1}}T_{c}{}_{b_{1}} 
-{1\over 2} R_{cb_{1}}T^{a_{1}c}
-R^{a_{1}}{}_{cdb_{1}} T^{cd}
\label{barDelta11}
\eea
\bea
-(\bar\D T)^{a_{1} a_{2}}{}_{b_{1}}=
\nabla^c\nabla_{c} \tilde{T}^{a_{1} a_{2}}{}_{b_{1}}
-R^{c[a_{1}}\tilde{T}_{c}{}^{a_{2}]}{}_{b_{1}}+
\frac{1}{2} 
R^{cda_{1}a_{2}}\tilde{T}_{cd}{}_{b_{1} }
\\ \nonumber \qquad -{1\over 2}R_{cb_{1}}\tilde{T}^{a_{1} a_{2}c}
-2R^{[a_{1}}{}_{cdb_{1}} \tilde{T}^{a_{2}]cd}
\label{barDelta21}
\eea
\bea
-(\bar\D T)^{a_{1} a_{2}}{}_{b_{1}b_{2}}=
\nabla^c\nabla_{c} \tilde{T}^{a_{1} a_{2}}{}_{b_{1}b_{2}}
-R^{c[a_{1}}\tilde{T}_{c}{}^{a_{2}]}{}_{b_{1} b_{2}}+
\frac{1}{2} 
R^{cda_{1}a_{2}}\tilde{T}_{cd}{}_{b_{1} b_{2}}
 \nonumber \\
-R_{c[b_{1}}\tilde{T}^{a_{1} a_{2}c}{}_{b_{2}]}+
\frac{1}{2} R_{cdb_{1}b_{2}} \tilde{T}^{a_{1} a_{2}cd}
-4R^{[a_{1}}{}_{cd[b_{1}} \tilde{T}^{a_{2}]cd}{}_{b_{2}]} 
\label{barDelta22}
\eea
\bea
-(\bar\D T)^{a_{1} a_{2}}{}_{b_{1}b_{2} b_{3}}=
\nabla^c\nabla_{c} \tilde{T}^{a_{1}a_{2}}{}_{b_{1}b_{2} b_{3}}
-R^{c[a_{1}}\tilde{T}_{c}{}^{a_{2}]}{}_{b_{1}b_{2} b_{3}}+
\frac{1}{2} 
R^{cda_{1}a_{2}}\tilde{T}_{cd}{}_{b_{1}b_{2} b_{3}}
 \nonumber \\
-\frac{3}{2}R_{c[b_{1}}\tilde{T}^{a_{1}a_{2}c}{}_{b_{2} b_{3}]}+
\frac{3}{2} R_{cd[b_{1}b_{2}} \tilde{T}^{a_{1}a_{2}cd}{}_{b_{3}]}
-6R^{[a_{1}}{}_{cd[b_{1}} \tilde{T}^{a_{2}]cd}{}_{b_{2} b_{3}]} 
\label{barDelta23}
\eea
\bea
-(\bar\D T)^{a_{1} a_{2}}{}_{b_{1}b_{2} b_{3}b_{4}}=
\nabla^c\nabla_{c} \tilde{T}^{a_{1}a_{2}}{}_{b_{1}b_{2} b_{3}b_{4}}
-R^{c[a_{1}}\tilde{T}_{c}{}^{a_{2}]}{}_{b_{1}b_{2} b_{3}b_{4}}+
\frac{1}{2} 
R^{cda_{1}a_{2}}\tilde{T}_{cd}{}_{b_{1}b_{2} b_{3}b_{4}}
 \nonumber \\
-2R_{c[b_{1}}\tilde{T}^{a_{1}a_{2}c}{}_{b_{2} b_{3}b_4]}+
3 R_{cd[b_{1}b_{2}} \tilde{T}^{a_{1}a_{2}cd}{}_{b_{3}b_{4}]}
-8R^{[a_{1}}{}_{cd[b_{1}} \tilde{T}^{a_{2}]cd}{}_{b_{2} b_{3}b_{4}]} 
\label{barDelta24}
\eea
 For the case of double $(q,p)$-forms, a 
straightforward computation provides the  commutation properties of the operators 
$d_{(i)}$ and $\delta_{(j)}$ for $i,j\in \{1,2\}$
\bea
\left([d_{(1)},d_{(2)}]T\right)^{a_{1}\dots a_{q+1}}{}_{b_{1}\dots b_{p+1}}=
\frac{(-1)^{p+q}}{2}(p+1)(q+1) \times\nonumber \\
\left(qR_{c[b_{p+1}}{}^{[a_{q}a_{q+1}}
\tilde{T}^{a_{1}\dots a_{q-1}]c}{}_{b_{1}\dots b_{p}]}-
pR^{c[a_{q+1}}{}_{[b_{p}b_{p+1}}
\tilde{T}^{a_{1}\dots a_{q}]}{}_{b_{1}\dots b_{p-1}]c}\right),\label{d12}
\eea
\bea
\left([d_{(1)},\delta_{(2)}]T\right)^{a_{1}\dots a_{q+1}}{}_{b_{1}\dots b_{p-1}}=
(-1)^{q} (q+1)\left(\frac{q}{2}R^c{}_{d}{}^{[a_{q}a_{q+1}}
\tilde{T}^{a_{1}\dots a_{q-1}]d}{}_{cb_{1}\dots b_{p-1}}\right.\nonumber \\
\left. +\frac{p-1}{2}R^{cd}{}_{[b_{1}}{}^{[a_{q+1}}
\tilde{T}^{a_{1}\dots a_{q}]}{}_{b_{2}\dots b_{p-1}]cd}+
R^{d[a_{q+1}}\tilde{T}^{a_{1}\dots a_{q}]}{}_{db_{1}\dots b_{p-1}}\right),
\label{ddelta}
\eea
\bea
\left([\delta_{(1)},\delta_{(2)}]T\right)^{a_{1}\dots a_{q-1}}{}_{b_{1}\dots b_{p-1}}
=\frac{p-1}{2}R^{ce}{}_{d[b_{1}}\tilde{T}^{da_{1}\dots a_{q-1}}{}_{b_{2}\dots 
b_{p-1}]ce} \nonumber \\-
\frac{q-1}{2}R_{ce}{}^{d[a_{1}}\tilde{T}^{a_{2}\dots a_{q-1}]ce}{}_{db_{1}\dots 
b_{p-1}}\label{delta12}
\eea
(Observe again, that in flat space these operators commute.)
In will be useful in subsequent calculations to note from (\ref{barDelta1T})  and (\ref{barDelta2T})
\bea
(\D_{(1)}T)^{a_{1}\dots a_{q}}{}_{b_{1}\dots b_{p}}-(\D_{(2)}T)^{a_{1}\dots a_{q}}{}_{b_{1}\dots b_{p}}= 
qR^{c[a_{1}}\tilde{T}_{c}{}^{a_{2}\dots a_{q}]}{}_{b_{1}\dots b_{p}}-
pR_{c[b_{1}}\tilde{T}^{a_{1}\dots a_{q}c}{}_{b_{2}\dots b_{p}]}\\ \nonumber
-\frac{q(q-1)}{2}
R^{cd[a_{1}a_{2}}\tilde{T}_{cd}{}^{a_{3}\dots a_{q}]}{}_{b_{1}\dots b_{p}}
+\frac{p(p-1)}{2} R_{cd[b_{1}b_{2}} 
\tilde{T}^{a_{1}\dots a_{q}cd}{}_{b_{3}\dots b_{p}]}
 \, .\label{barDelta1-barDelta2}
\eea

\section{Double $(p,p)$-forms: Curvature tensors}
The transpose ${}^t T$ of $T$ is of special relevance for the special
case of the double $(q,p)$-forms in which both blocks have the same
number of indices, say $p$ ---these double $(p,p)$-forms include, of
course, the important case of curvature tensors.  In this case the two
blocks of $T$ can be interchanged, and hence we will say that a double
$(p,p)$-form is {\em symmetric} if $T={}^tT$, and {\it antisymmetric}
if $T=-{}^tT$ (in this case only for $p>1$).  Of course, for $p>1$ any
double $(p,p)$-form can be decomposed uniquely into a symmetric and an
antisymmetric one (if $p=1$ the decomposition gives a 2-form and a
double symmetric $(1,1)$-form).  Hence, without loss of generality, in
what follows we will only consider these two cases,
$$
T = \pm {}^tT  \qquad \hbox{or} \qquad\tilde{T}^{a_{1}\dots a_{p}}{}_{b_{1}\dots b_{p}}=
\pm \tilde{T}_{b_{1}\dots b_{p}}{}^{a_{1}\dots a_{p}}\, .
$$
Then it follows trivially from (\ref{trans-ddel}) that 
\be
d_{(2)}T=\pm\, {}^t (d_{(1)}T) ,\hspace{1cm} \delta_{(2)}T =\pm\,  {}^t (\delta_{(1)}T)
\label{trans-d}
\ee
%% 
 % $$
 % (d_{(2)}T)^{a_{1}\dots a_{p}}{}_{b_{1}\dots b_{p+1}}=\pm (d_{(1)}T)_{b_{1}\dots 
 % b_{p+1}}{}^{a_{1}\dots a_{p}}, \,\,\,
 % (\delta_{(2)}T)^{a_{1}\dots a_{p}}{}_{b_{1}\dots b_{p-1}}=\pm
 % (\delta_{(1)}T)_{b_{1}\dots b_{p-1}}{}^{a_{1}\dots a_{p}}
 % $$
 %%
and therefore $\D_{(2)}T=\pm\, {}^t(\D_{(1)}T)$ so that
$$
\bar\D T =\frac{1}{2}\left(\D_{(2)}T\pm {}^t(\D_{(2)}T)\right).
$$
Since the associated superpotential $\stackrel{o}{\mbox{$T$}}$ will have the
same symmetry properties as $T$, it follows that the four potentials
$Y_{(1)}$, $Y_{(2)}$, $Z_{(1)}$ and $Z_{(2)}$ defined in the previous
section satisfy
\bea
Y_{(2)}&=&\pm ^t Y_{(1)}\equiv (-1)^{p+1} Y,\label{Y}\\
Z_{(2)}&=&\pm ^t Z_{(1)}\equiv  (-1)^{p-1}Z \, ,\label{Z}
\eea
where furthermore the completely antisymmetric part ${\cal A}[Y_{(2)}]$ 
of $Y_{(2)}$ (or 
$Y$) may vanish identically due to the (anti)symmetry of the associated 
superpotential:
\bea
{\cal A}[Y_{(2)}]=0 \,\, \mbox{for}\,\,\left\{
\begin{array}{lr}
      T={}^tT  & p\,\, \mbox{odd}, \\ 
    T=-{}^t T & p\,\,  \mbox{even}.
\end{array} \right.
\label{Ysym}
\eea
The index notation for this restriction is simply
\bea
Y_{[a_1\ldots a_pb_1\ldots b_{p+1}]}=0 \,\, \mbox{for}\,\,\left\{
\begin{array}{lr}
      T={}^tT  & p\,\, \mbox{odd}, \\ 
    T=-{}^t T & p\,\,  \mbox{even}.
\end{array} \right.
\label{Ysymi}
\eea

\begin{theorem}
\label{theo:Y+Z}
Given any tensor $T$ with the structure of a
double (anti)symmetric $(p,p)$-form 
there  always exists an {\em associated} local superpotential 
${\stackrel{o}{\mbox{$T$}}}$ such that
$ \bar\D \stackrel{o}{\mbox{$T$}} =T$; 
furthermore there always exist a {\em pair} of local
potentials $Y_{(2)}, Z_{(2)}$ satisfying (\ref{Ysym}) such that
\be
T=\frac{1}{2}\left[\delta_{(2)}Y_{(2)}\pm{}^t(\delta_{(2)}Y_{(2)})+
d_{(2)}Z_{(2)}\pm{}^t(d_{(2)}Z_{(2)})\right],\label{4pota}
\ee
where the double $(p,p+1)$-form $Y_{(2)}=d_{(2)}\stackrel{o}{\mbox{$T$}}$, 
and the double $(p,p-1)$-form $Z_{(2)}=\delta_{(2)}\stackrel{o}{\mbox{$T$}}$. 

Using the notation introduced in
(\ref{Y}) and (\ref{Z}), (\ref{4pota}) can be written in index notation as
\bea
\tilde{T}^{a_{1}\dots a_{p}}{}_{b_{1}\dots b_{p}}&=&\frac{1}{2}\left(
Y^{a_{1}\dots a_{p}}{}_{b_{1}\dots b_{p}c}{}^{;c} \pm
Y_{b_{1}\dots b_{p}}{}^{a_{1}\dots a_{p}c}{}_{;c}\right. \nonumber \\
&& + \left. pZ^{a_{1}\dots a_{p}}{}_{[b_{1}\dots b_{p-1};b_{p}]}\pm
pZ_{b_{1}\dots b_{p}}{}^{[a_{1}\dots a_{p-1};a_{p}]}\right)
\label{T=Y+Z}
\eea
where  the potential $Y$ satisfies (\ref{Ysymi}), and the 
potentials themselves can be given in terms of 
the double (anti)symmetric associated superpotential 
$\stackrel{o}{\mbox{$T$}}$ by 
$$
Y^{a_{1}\dots a_{p}}{}_{b_{1}\dots b_{p+1}}=-(p+1) 
{\stackrel{o}{\mbox{$T$}}}{}^{a_{1}\dots a_{p}}{}_{[b_{1}\dots b_{p};b_{p+1}]}, 
\qquad 
Z^{a_{1}\dots a_{p}}{}_{b_{1}\dots b_{p-1}}=
-{\stackrel{o}{\mbox{$T$}}}{}^{a_{1}\dots a_{p}}{}_{c b_{1}\dots  
b_{p-1}}{}^{;c}.
$$     
\end{theorem}

Since the operations of transpose and trace 
commute, then if $T$ is (anti)symmetric, so is
$\hbox{tr}(T)$.  Suppose now that, in addition to the (anti)symmetry
between blocks, the double $(p,p)$-form $T$ is {\em traceless}, i.e.,
\be
\hbox{tr}(T)=0\, .\label{Ttraceless}
\ee
The associated superpotential ${\stackrel{o}{T}}$ will also be
{traceless} with
\be
\hbox{tr}(\stackrel{o}{\mbox{$T$}})=0\, .\label{Totraceless}
\ee
But then the potentials $Y_{(2)}$ and $Z_{(2)}$ of (\ref{4pota}) 
{\em are not independent} and 
one can easily check using (\ref{Totraceless}) that $Z_{(2)}$ is essentially the trace of  the transpose of 
$Y_{(2)}$:
$$
Z_{(2)}=\mp \ \hbox{tr}({}^tY_{(2)})
$$
or with indices and the notation of (\ref{Y}-\ref{Z})
\be
Z^{a_{1}\dots a_{p}}{}_{b_{1}\dots b_{p-1}}=
\mp \, Y_{cb_{1}\dots b_{p-1}}{}^{ca_{1}\dots a_{p}} \, . \label{trY=Z}
\ee
Furthermore
$$
\hbox{tr}(\hbox{tr}(Y_{(2)}))=0
$$
which becomes in index notation
\be
Y_{cdb_{1}\dots b_{p-2}}{}^{cda_{1}\dots a_{p-1}}=0\, . \label{trY=0}
\ee
Therefore, we have
\begin{theorem}
\label{theo:Y}
Given any tensor $T$ with the structure of a double (anti)symmetric 
traceless $(p,p)$-form there always exists an {\em associated}
local superpotential ${\stackrel{o}{\mbox{$T$}}}$ such that
$\bar\D \stackrel{o}{\mbox{$T$}} =T$; 
furthermore, there always exists a double $(p,p+1)$-form local
potential $Y_{(2)}$ satisfying (\ref{Ysym}) and 
$\mbox{\rm{tr}}(\mbox{\rm{tr}}(Y_{(2)}))=0$ such that 
\be
T=\frac{1}{2}\left(\delta_{(2)}Y_{(2)}\pm {}^t(\delta_{(2)}Y_{(2)})
-d_{(1)}\mbox{\rm{tr}}(Y_{(2)}) \mp{}^t(d_{(1)}\mbox{\rm{tr}}(Y_{(2)}))\right),
\label{T=Ya}
\ee
where the double $(p,p+1)$-form $Y_{(2)}=d_{(2)}\stackrel{o}{\mbox{$T$}}$.

With the notation introduced in (\ref{Y}), the index version of
(\ref{T=Ya}) can be written as
\bea
\tilde T^{a_{1}\dots a_{p}}{}_{b_{1}\dots b_{p}}&=&\frac{1}{2}\left(
Y^{a_{1}\dots a_{p}}{}_{b_{1}\dots b_{p}c}{}^{;c} \pm
Y_{b_{1}\dots b_{p}}{}^{a_{1}\dots a_{p}c}{}_{;c}\right. \nonumber \\
&&\qquad  -
pY^{c[a_{1}\dots a_{p-1}}{}_{cb_{1}\dots b_{p}}{}^{;a_{p}]} \mp 
\left. pY_{c[b_{1}\dots b_{p-1}}{}^{ca_{1}\dots a_{p}}{}_{;b_{p}]}\right)
\label{T=Y}
\eea
where the potential $Y$ satisfies (\ref{Ysymi}) and (\ref{trY=0}), and
is given in terms of the double 
(anti)symmetric traceless associated
superpotential ${\stackrel{o}{\mbox{$T$}}}$ by
\be
Y^{a_{1}\dots a_{p}}{}_{b_{1}\dots b_{p+1}}\equiv -(p+1)
{\stackrel{o}{T}}{} ^{a_{1}\dots a_{p}}{}_{[b_{1}\dots 
b_{p};b_{p+1}]}\, .\label{Y=d2t}
\ee 
\end{theorem}

In \cite{E3} it was shown in $2m$ dimensions that a traceless
symmetric $(m,m)$-form has a traceless $(m,m-1)$-form potential.  That
result is a special case of the much more general result in this
theorem, where the double dual of the $(m,m-1)$-form potential in
\cite{E3} is precisely the $(m,m+1)$-form potential $Y$ in this
theorem.

For future reference we give explicitly the index version of the
potential structure for a symmetric $2$-tensor
\bea
{T}_{a}{}_{b}&=&\frac{1}{2}\left(
Y_{a}{}_{bc}{}^{;c} +
Y_{b}{}_{ac}{}^{;c}\right.  + \left. Z_{a}{}_{;b}+
Z_{b}{}_{;a}\right).
\label{T2=Y+Z}
\eea
and for a {\it traceless} symmetric $2$-tensor
\bea
{T}_{a}{}_{b}&=&\frac{1}{2}\left(
Y_{a}{}_{bc}{}^{;c} +
Y_{b}{}_{ac}{}^{;c}  -Y^c{}_{ca;b}-
Y^c{}_{cb;a}\right).
\label{T2=Y}
\eea
where (in both cases) the double $(1,2)$-form $Y^a{}_{bc}$  
satisfies $Y_{[abc]}=0$.
 
There is a well known decomposition for symmetric $2$-tensors in {\it
three} dimensional spaces, but it is restricred to {\it proper Riemannian space} \cite{De, Yo},  and so differs from this one.  Of
course, the antisymmetric $2$-tensor is just a single $2$-form whose
decomposition (\ref{r2eq}) is easily seen to agree with the
decomposition which follows from Theorem \ref{theo:Y+Z}.
 
\subsection{Application to general curvature tensors}
Let us apply the above results to the case of {\em Riemann candidates}, 
that is, tensors with the algebraic properties of a Riemann curvature 
tensor. Let ${\cal R}_{abcd}$ be any such Riemann candidate, that is to say
\be
{\cal R}_{abcd}={\cal R}_{[ab][cd]}, \hspace{1cm} {\cal R}_{a[bcd]}=0 \,\,\,
(\Longrightarrow \,\,\, {\cal R}_{abcd}={\cal R}_{cdab}),\label{F-prop}
\ee
so that ${\cal R}_{abcd}$ is in particular a {\em symmetric} double (2,2)-form.

Let ${\stackrel{o}{{\cal R}}}{}_{abcd}$  
be the associated local superpotential for
${\cal R}_{abcd}$ by
\be                                                                    
(\bar\D {\stackrel{o}{\mbox{${\cal R}$}}}){}_{abcd}=
{\cal R}_{abcd} \label{boxF=F}                   
\ee   
and because of the properties of $\bar\Delta$, 
${\stackrel{o}{\mbox{${\cal R}$}}}$ is also a Riemann candidate;  
the pair of potentials defined in Theorem \ref{theo:Y+Z} become 
now
\be
Y_{abcde}=Y_{[ab][cde]}\equiv
-3 {\stackrel{o}{\mbox{${\cal R}$}}}{}_{ab[cd;e]}, \hspace{1cm}
Z_{abc}= Z_{[ab]c}\equiv -
{\stackrel{o}{\mbox{${\cal R}$}}}{}_{abdc}{}^{;d} \, .\label{Y=F}
\ee
It is important to realize, since ${\stackrel{o}{\mbox{${\cal R}$}}}$ 
is a Riemann candidate satisfying (\ref{F-prop}), that these potentials satisfy 
the additional symmetries
\be
Y_{a[bcde]}=0, \hspace{1cm} Z_{[abc]}=0  , \label{Ypty1}
\ee
the first of which implies the following useful properties 
\be
Y^e{}_{[bcd]e}=0, \hspace{5mm} Y_{[abcd]e}=0, \hspace{5mm}
Y_{abcde}=3 Y_{[cde]ab}=3Y_{a[cde]b}, \hspace{5mm} 
Y_{a[bc]de}=-Y_{a[de]bc} .\label{Ypty2}
\ee

Given that (\ref{boxF=F}) always has local solutions 
${\stackrel{o}{\mbox{${\cal R}$}}}$ for any 
given ${\cal R}$ we have,
\begin{theorem}
\label{theoremR}
Any Riemann candidate tensor  ${\cal R}^{ab}{}_{cd}$ has a pair
of local potentials given by a double $(2,3)$-form $Y^{ab}{}_{cde}$
and a double (2,1)-form $Z^{ab}{}_{c}$ with the properties
(\ref{Ypty1}) such that
\be                                                                    
{\cal R}_{abcd}=\frac{1}{2}\left(Y_{abcde}{}^{;e}+Y_{cdabe}{}^{;e}+       
2Z_{ab[c;d]}+2Z_{cd[a;b]}\right)\, .\label{G=Y+Pi}                                     
\ee   
The potentials themselves can be given in terms of the associated
(Riemann candidate) local superpotential 
${\stackrel{o}{{\cal R}}}{}_{abcd}$ by (\ref{Y=F}).
\end{theorem}

Of course, any Riemann candidate can be decomposed in terms of its 
trace (its `Ricci tensor') and its traceless part (its `Weyl 
tensor'). The trace is a double symmetric (1,1)-form, and therefore 
the general Theorem \ref{theo:Y+Z} applies, and in particular (\ref{T2=Y+Z}). After a little rearranging, it can be 
seen that the potentials for this Ricci part are essentially the 
traces of the pair of potentials  defined in (\ref{Y=F}) and the formula 
relating them is the trace of (\ref{G=Y+Pi}).
\begin{coro}
\label{corRic}
Any Ricci candidate tensor  ${\cal R}_{ab}$ has a pair
of local potentials given by a double $(1,2)$-form $Y^{c}{}_{ab}$ 
with the property $Y_{[abc]}=0$ and a double (1,0)-form $Z_{a}$
such that
\be                                                                    
{\cal R}_{ab}=\frac{1}{2}\left(Y_{abe}{}^{;e}+ Y_{bae}{}^{;e}+       
Z_{a;b}+Z_{b;a}\right)\, .\label{Ric=Y+Z}                                     
\ee   
The relation with the Riemann candidate local potentials of (\ref{G=Y+Pi}) is
$$Y^{c}{}_{ab} = Y^{ec}{}_{eab} -Z_{ab}{}^c, \hspace{1cm}
Z_{a} =Z_{ea}{}^e\, .
$$
\end{coro}
The traceless part of the Ricci candidate tensor (essentially any
traceless symmetric $2$-tensor) has an even simpler structure
requiring only one potential, as shown in (\ref{T2=Y}):
\begin{coro}
\label{corRict}
Any traceless Ricci candidate tensor 
$\hat{\cal R}_{ab} ={\cal R}_{ab}-\frac{1}{n}g_{ab}{\cal R}$, where ${\cal R}={\cal R}^c{}_{c}
$, 
has a local potential given by a double $(1,2)$-form $\hat Y^{c}{}_{ab}$
 with the property
$\hat Y_{[abc]}=0$ such that
\be                                                                    
\hat {\cal R}_{ab}=\frac{1}{2}\left(\hat Y_{abe}{}^{;e}+\hat Y_{bae}{}^{;e}-      
\hat Y^{c}{}_{ca;b}-\hat Y^{c}{}_{cb;a}\right)\, .\label{Rict=Y}                                     
\ee 
The relation with the Ricci candidate local potentials of (\ref{Ric=Y+Z}) is
$$\hat Y^{c}{}_{ab} = Y^{c}{}_{ab} +\frac{2}{n}g_{c[a}Y^e{}_{b]e} +\frac{2}{n}g_{c[a} Z_{b]}\, .
$$  
\end{coro}

More interesting is the traceless part of a Riemann candidate, 
which we call a {\em Weyl candidate}, that is,  
a double (2,2)-form ${\cal C}_{abcd}={\cal C}_{[ab]cd}={\cal C}_{ab[cd]}$ 
with the algebraic properties of the Weyl conformal curvature tensor:
\be
{\cal C}^a{}_{bca}=0 , \hspace{5mm} {\cal C}_{a[bcd]}=0,  \,\,\,
(\Longrightarrow \,\,\, {\cal C}_{abcd}={\cal C}_{cdab}),\label{v-prop}
\ee
so that ${\cal C}_{abcd}$ is in particular
a {\em traceless and symmetric} double (2,2)-form. 
By considering the associated superpotential 
$\stackrel{o}{\mbox{${\cal C}$}}$ for ${\cal C}$ given by 
${\cal C} = \bar\D {\stackrel{o}{\mbox{${\cal C}$}}}$,
we can then deduce from  Theorem 
\ref{theo:Y} and Theorem \ref{theoremR}, 
\begin{theorem}
\label{theoremW}
Any Weyl candidate tensor field
${\cal C}_{abcd}$ has a double $(2,3)$-form local potential $P_{abcde}$ 
with the properties 
\be
P_{a[bcde]}=0, \hspace{1cm} P^{ab}{}_{abc}=0 \label{Ppty1}
\ee
such that
\be 
{\cal C}^{ab}{}_{cd} = 
\frac{1}{2}\left(P^{ab}{}_{cde}{}^{;e}+P_{cd}{}^{abe}{}_{;e} -2 
P_e{}_{[c}{}^{abe}{}_{;d]}-2P^e{}^{[a}{}_{cde}{}^{;b]} \right)\, . \label{W=P}
\ee
The potential itself can be given in terms of an associated local 
(Weyl candidate)
superpotential ${\stackrel{o}{\mbox{${\cal C}$}}}{}^{ab}{}_{cd}$ by 
\be
P^{ab}{}_{cde}=P^{[ab]}{}_{[cde]}\equiv - 3 
{\stackrel{o}{\mbox{${\cal C}$}}}{}^{ab}{}_{[cd;e]} \label{P=W}\ .
\ee
\end{theorem}

This result was originally obtained, by a more direct route, in
\cite{ES}.  The operator $\bar\D$ which we have used in this theorem,
specialising via (\ref{barDelta22}), is easily seen to coincide with
the operator which we constructed for the superpotential of the Weyl
candidate tensor in \cite{ES}.
As noted there, in four dimensions our new
(2,3)-form potential $P^{ab}{}_{cde}$ coincides with the double dual
of the classical Lanczos (2,1)-form potential $H^{ab}{}_c$ for the Weyl tensor
\cite{La}.

Immediate consequences from the first of (\ref{Ppty1}) are the following 
useful properties, reminiscent of (\ref{Ypty2})
\be
P^e{}_{[bcd]e}=0, \hspace{5mm} P_{[abcd]e}=0, \hspace{5mm}
P^{ab}{}_{cde}=3 P_{[cde]}{}^{ab}=3 P^{[a}{}_{[cde]}{}^{b]}, \hspace{5mm} 
P_{a[bc]de}=-P_{a[de]bc}\, . \label{Ppty2}
\ee

The basic equation (\ref{W=P}) can be written in several alternative 
but equivalent forms by using properties (\ref{Ppty1}) and (\ref{Ppty2}). 
Other possibilities are 
\be 
{\cal C}^{ab}{}_{cd} = P^{ab}{}_{cde}{}^{;e} + P^{e[ab]}{}_{cd;e} - 
P_e{}_{[c}{}^{abe}{}_{;d]}-P^e{}^{[a}{}_{cde}{}^{;b]}  \label{W=P'}
\ee
and
\be
{\cal C}^{ab}{}_{cd} = P^{ab}{}_{cde}{}^{;e} + P^{e[ab]}{}_{cd;e}
-P^e{}^{[a}{}_{cde}{}^{;b]}-2P^{e[ab]}{}_{e[c;d]}   \label{W=P''}.
\ee

\section{Laplace-type equation for the Riemann tensor}
\label{sec:Riemann}
The intrinsic, as well as quickest and simplest, way to obtain a 
Laplace-like equation (a wave equation in Lorentzian signature) for 
the electromagnetic field ${F}$ is as follows. The 2-form ${F}$ 
satisfies the Maxwell equations
\be
d{F}=0 , \hspace{1cm} \delta{F} = J \label{maxwell}
\ee
where $J$ is the electric current 1-form. By using the de Rham 
operator the Laplace-like equation follows from these equations 
strightforwardly 
\be
\Delta { F}= d J \, .\label{DF}
\ee
As is known, in the absence of electric charges and currents ($J=0$) one has 
the simpler equation
$$
\mbox{ $J=0$} \hspace{3mm} \Longrightarrow \hspace{2mm} \D 
{F}=0 \, .
$$
The index version of (\ref{DF}) is the well-known equation (``wave'' equation 
in Lorentzian signature) 
$$
\nabla^c\nabla_{c} F^{ab}+R^{efab}F_{ef}-R^{ae}F_{e}{}^b + 
R^{be}F_{e}{}^a=J^{[a;b]}
$$
or equivalently
$$
\nabla^c\nabla_{c} F^{ab}-2R^{eafb}F_{ef}-R^{ae}F_{e}{}^b + 
R^{be}F_{e}{}^a=J^{[a;b]}\, .
$$

As has been always the case, the electromagnetic field may serve as 
source of inspiration for the gravitational field, and thereby for the 
curvature tensor of any pseudo-Riemannian manifold. Thus, there have 
been some attempts to formulate the analogue of (\ref{DF})
for the Riemann tensor with
a generalised de Rham operator, using  the Bianchi identities as 
the analogue to the Maxwell equations (\ref{maxwell}). 
We will use $\Re $ to represent the double symmetric 
(2,2)-form defined by the Riemann tensor
$R^{ab}{}_{cd}$. 
%% 
 % as a $2$-tensor valued $2$-form with the last pair of
 % indices treated as the form indices, and the first pair as the tensor
 % indices.  
 %%
Then, the Bianchi identities $R_{ab[cd;e]}=0$ can be written
in a concise form using the notation of this paper as
\be
d_{(2)}\Re=0\, . \label{bianchi}
\ee
%% 
 % We emphasise that from our point of view $\Re$ can also be viewed as a
 % $2$-tensor-valued $2$-form with the first pair of indices in
 % $R^{ab}{}_{cd}$ as form indices, and the last pair as ordinary tensor
 % indices;
 %%
Hence, using (\ref{trd}) to
take the trace of (\ref{bianchi}) we readily obtain
\be
\delta_{(1)}\Re = - d_{(2)}\Re ic \ (= d_{(2)}\hbox{tr}(\Re )) \label{tr-bianchi}
\ee
where $\Re ic = \hbox{tr}(\Re )$ denotes the double symmetric (1,1)-form 
representing the Ricci tensor $R_{ab}$. 
(Alternatively (\ref{tr-bianchi}) could be obtained by translating
from the contracted index version of the Bianchi identities
$R_{abcd}{}^{;a}=-2R_{b[c;d]}$.)
%% 
 % ; but we prefer working within the
 % tensor-valued form formalism which avoids possible ambiguities.)
 %%

The symmetry  of the Riemann 
tensor ${}^t\Re =\Re $ implies that (\ref{bianchi}) is actually 
equivalent, via (\ref{trans-d}), to 
\be
d_{(1)}\Re =0 
\label{bianchi2}
\ee
and by taking the trace we obtain  \be
\delta_{(2)}\Re = - d_{(1)} \Re ic \ \ 
\label{tr-bianchi2}
\ee
(or equivalently using the symmetry of the Riemann tensor in  (\ref{tr-bianchi})).

By following a 
procedure similar to that leading to (\ref{maxwell}) 
exploiting the de Rham operator for tensor valued 2-forms 
\cite{BCJR,BCJF,MTW,BF}, which in our notation is simply $\D_{(2)}$ (or 
$\D_{(1)}$), by using 
(\ref{bianchi},\ref{tr-bianchi2}), one easily derives
\be
\D_{(2)}\Re = -d_{(2)}d_{(1)}\Re ic \, .\label{DR=R}
\ee
For those spaces with $\nabla_{[a}R_{b]c}=0$,  
\be
\mbox{ $d_{(1)}\Re ic=0$} \hspace{3mm} \Longrightarrow \hspace{2mm} \D_{(2)} 
\Re=0 \, . \label{harmonic2}
\ee
This Laplace-like equation for the Riemann tensor appears to have been
written down first of all by Penrose \cite{Pe}, motivated by spinors,
and later in \cite{MTW}; in both cases, only four dimensional 
Ricci-flat spaces were considered.  Subsequently the general version (\ref{DR=R})
has been
presented in \cite{Bel,Ry,MT,BCJF}, as a direct result of
differentiating the Bianchi equations. The typical index version of 
(\ref{DR=R}) reads
\be
\label{box2R=RR}
\nabla^e\nabla_{e} {  R}^{ab}{}_{cd}+
4{R}^{[a}{}_{ef[c}{}{R}^{b]ef}{}_{d]}
+{R}^{abef}{ R}_{ef}{}_{cd}-2{R}^{ab}{}_{e[c}R^e{}_{d]}=
4{  R}^{[a}{}_{[c;}{}^{b]}{}_{d]} \, .\label{boxR1}
\ee
(One should note different sign conventions, and some possible sign
inconsistencies in the literature; see \cite{AE1} for a summary.)

As explained earlier in Section \ref{sec:r-fold}, 
the trace and index symmetry properties across the two
pairs of indices do not commute with the operator $\D_{(2)}$; and so
in particular, 
\bean
(i)\hspace{2mm} (\D_{(2)}{\Re})_{a[bcd]} \ne 0, \hspace{2cm}
(ii)\hspace{2mm} \hbox{tr}( \D_{(2)}{\Re}) \neq \D_{(2)}\hbox{tr}({\Re})
=\D_{(2)} {\Re ic}
%(\D_{(2)}\bm{R})^{ab}{}_{cb} \ne(\D_{(2)}R)^{a}{}_{c}
\eean
as is easily confirmed.Ê However, this operator $\D_{(2)}$ has been proposed 
in \cite {MTW,BF,BCJR,BCJF} as the natural analogue for the
Riemann tensor of the usual de Rham operator $\D$ for the electromagnetic
field tensor; we would argue, precisely because of these deficiencies,  
that it is not the most appropriate operator. 
It should be understood that in these references only
vacuum spaces are being explicitly considered, and in such spaces the
commutativity properties $(ii)$ which fail in general are satisfied trivially
because of the absence of Ricci tensor terms. 

Of course, we could have done the whole calculation using the 
corresponding transposes, and we would easily arrive at 
\be
\D_{(1)}\Re = -d_{(1)}d_{(2)}\Re ic \label{D2R=R}\ 
\ee
which is strictly equivalent to (\ref{DR=R}), 
due to the symmetry of the Riemann tensor pairs.

A similar argument could be given concerning the information
obtainable using the weighted operator $\bar\D$.  Despite this fact,
in general we believe that the operators $\D_{(1)}$ and $\D_{(2)}$ are
potentially confusing and not really particularly suitable; this is
especially true for general tensors with no symmetry between blocks,
but also for the Riemann tensor itself.  We firmly believe that the
proper operator to be used here is $\bar\D$.  From (\ref{DR=R},
\ref{D2R=R}) it follows immediately that
\be
\bar\Delta \Re=  Ê-Ê\frac{1}{2}(d_{(1)}d_{(2)} +d_{(2)}d_{(1)})\Re ic\, .
\label{barDRf}
\ee
This is a more useful version of a Laplace-like equation for the
Riemann tensor, and it does not have the disadvantages of
(\ref{DR=R}).  It is easy to check, for instance, that the two
examples considered above lead immediately to trivial identities:
\bean
(i)\hspace{2mm} (\bar\Delta {\Re})_{a[bcd]} = 0, \hspace{2cm}
(ii)\hspace{2mm} \hbox{tr}(\bar\Delta {\Re}) =\bar\D 
\hbox{tr}({\Re}) =\bar\D {\Re ic} \, .
\eean
It should be remarked that the second equation here implies, in fact,
that there is no non-trivial expression for $\bar\D \Re ic$ derivable
from (\ref{barDRf}).  If one wishes to find a Laplace-like equation
for the Ricci tensor, the best route is to use (\ref{trd},
\ref{trdelta}) to take the trace of (\ref{tr-bianchi}) to 
get\footnote{This formula (\ref{tr-tr-bianchi}) is an illustrative 
example of the possibility mentioned on a number of occasions in this paper that some 
multiple forms must be considered as having extra blocks with no indices, when combined with other such 
forms. In this case, despite 
$R$ being a scalar, one has to take its ``$d_{(2)}$''-derivative. This 
means that $R$ in this equation acts as a `double (0,0)-form'. The 
general rule for these situations is to define the form-structure 
number of an equation as the maximum one for the different terms 
involved, and then every term in such equation has to be considered 
accordingly. See also 
Eqs.(\ref{D1Ric}-\ref{barDRscalar}), (\ref{Ric=Y+Z})
for other examples.}
\be 
\delta_{(1)}\Re ic = - \frac{1}{2} d_{(2)} R \label{tr-tr-bianchi}
\ee
where $R=\hbox{tr}({\Re})=R^c{}_c$ is the scalar curvature;
from which, together with (\ref{tr-bianchi2}) it follows that
\be
\D_{(1)}\Re ic = -\delta_{(1)}\delta_{(2)} \Re - \frac{1}{2} d_{(1)}  d_{(2)} R 
.\label{D1Ric}
\ee
For the same reasons as above we prefer to combine it with its $\D_{(2)}$
counterpart
\be
\D_{(2)} \Re ic = -\delta_{(2)} \delta_{(1)} \Re ic - \frac{1}{2} d_{(2)} d_{(1)}  R 
\label{D2Ric}
\ee
to obtain
\be
\bar\D\Re ic =-\frac{1}{2}(\delta_{(1)}\delta_{(2)}+\delta_{(2)}\delta_{(1)}){\Re}-\frac{1}{4}(d_{(1)} d_{(2)} + d_{(2)}  d_{(1)} ) R \, .
\label{barDRic}
\ee
A final contraction here can be obtained by remembering that $\bar \D$
commutes with contraction, and using (\ref{trdelta}) on the first term
on the right hand side; on the otherhand, we cannot use (\ref{trd}) on
the second term since $R$ is a scalar ---acting here as a double $(0,0)$-form. 
However, in index notation $[(d_{(1)} d_{(2)} + d_{(2)} d_{(1)} ) R]_{ab} = R_{;ab}+R_{;ba}$ 
and hence its trace is $2\bar\D R $.  
Hence we obtain the usual Laplace-like equation
for the scalar curvature $R$
\be
\frac{1}{2}\bar\D R =- \frac{1}{2}(\delta_{(1)}\delta_{(2)}+\delta_{(2)}\delta_{(1)})  \Re ic =
- \delta_{(1)}\delta_{(2)}  \Re ic \, ,\label{barDRscalar}
\ee
where the last equality is due to the symmetry of the Ricci tensor.

The index version of (\ref{barDRf}) was written in \cite{DW},
\cite{AE1} and reads
\begin{eqnarray}
\label{boxR=RR} 
&\nabla^e\nabla_e {  R}^{ab}{}_{cd}+4{R}^{[a}{}_{ef[c}{}{R}^{b]ef}{}_{d]}+
{  R}^{abef}{R}_{ef}{}_{cd}+{  R}^{ab}{}_{e[c}R^e{}_{d]}+
{  R}_{cd}{}^{e[a}R_e{}^{b]} \\ \nonumber& =
2{  R}^{[a}{}_{[c}{}^{;b]}{}_{d]}+2{  R}_{[c}{}^{[a}{}_{;d]}{}^{b]} 
\label{boxR2}
\end{eqnarray} 
while the index versions of (\ref{barDRic}) and (\ref{barDRscalar})
are the familiar equations:
$$
\nabla^e\nabla_e R_{ab} -\frac{1}{2} \nabla_a\nabla_b R-R_a{}^{c}R_{bc}-R_{aefb}R^{ef}=\nabla_e\nabla_fR^{e}{}_{(a}{}^f{}_{b)},
$$
$$
 \frac{1}{2} \nabla^e\nabla_e R =\nabla_e\nabla_fR^{ef}\ .
$$
We believe that the equation (\ref {barDRf})
is the natural Riemann tensor analogue of (\ref{DF}), and of course we also have
\be
\mbox{$d_{(1)}\Re ic =0$} \hspace{3mm} \Longrightarrow \hspace{2mm} \bar\D
{\Re}=0 \, .\label{harmonicbar}
\ee
Let us comment on the implications of the different definitions of
`harmonic tensors' for the Riemann tensor.  From (\ref{harmonic2}), 
(\ref{harmonicbar}) we see 
that the Riemann tensor is fully harmonic (ergo harmonic) in 
spaces with $d_{(1)}{\Re ic} =0$ ($\Longleftrightarrow d_{(2)}{\Re 
ic}=0$), including manifolds with a parallel Ricci tensor and its 
particular cases of Ricci-flat and Einstein spaces. As noted above, 
it has
sometimes been the case that the operator $\D_{(2)}$ has been
considered as the extension to arbitrary tensors of the de Rham
operator for single forms; and hence for the special case of the
Riemann tensor the condition $\D_{(2)} \Re = 0 $ has been identified
as its harmonic condition \cite{BF, BCJF}. However, 
we want to stress the following two important remarks, which 
highlight the properties of $\bar\D$ and support our claim that this 
is the right operator, giving the proper harmonic condition:
\begin{itemize}
\item The class of pseudo-Riemannian manifolds with a harmonic Riemann 
tensor is {\em larger} than that containing a fully harmonic Riemann 
tensor. This is clear from (\ref{DR=R}), (\ref{D2R=R}) and 
(\ref{barDRf}), as we only need the condition 
$(d_{(1)}d_{(2)}+d_{(2)}d_{(1)}){\Re ic}=0$ for the former case. (Of course, 
for compact without boundary proper Riemannian manifolds the fully 
harmonic and harmonic cases are equivalent, as follows from Theorem 
\ref{harmonic=full}.)
\item We also note that the most common Laplacian used for
arbitrary tensors, and thought to be mathematically appropriate,
has been the Lichnerowicz Laplacian $\D_L$. We
have pointed out, in general, that tensors harmonic in the sense of 
$\D_{L}T=0$, do not coincide with (fully) harmonic tensors; 
in particular, as follows from (\ref{relation}),
{\em the Riemann tensor of pseudo-Riemannian manifolds 
(compact or not), are not harmonic by using $\D_{L}$, not even for 
Ricci-flat cases}. Probably this was 
the reason that some authors 
chose equation (\ref{DR=R}) (or equivalently (\ref{D2R=R})) to define 
the Laplace-like or wave equations for the Riemann tensor, letting 
aside the probably `more natural' operator $\D_{L}$.
\end{itemize}
We believe that these two remarks taken together with other 
advantages mentioned before demonstrate that the 
weighted operator $\bar\D$ is better suited than either of 
$\D_{(i)},\D_{L}$, and probably the right 
choice for defining harmonic tensors and studying the derived 
implications.

\

Finally, we want to remark about yet another good property of $\bar\D$
(shared in this case by $\D_{(1)}$ and $\D_{(2)}$): for any double
(2,2)-form ${\cal R}$ satisfying the first Bianchi identity (that is,
a Riemann candidate), its `Weyl part' is uniquely defined as usual by
taking out the traces.  We can denote by $W\{{\cal R}\}$ the Weyl part
of a Riemann candidate.  For example, $W\{{\Re}\}= {C}$ where $C$
denotes the Weyl tensor.  As proved first in \cite{Bel}, the operator
$W$ commutes with $\D_{(2)}$, and therefore with $\bar\D$.  Hence, this
provides the straightest route to the correct Laplace-like equation
for the Weyl tensor, because
$$
W\{\bar\D{\Re}\}=\bar\D(W\{{\Re}\})=\bar\D C
$$
so that, on using (\ref{barDRf}) we obtain
$$
\bar\D C= -Ê\frac{1}{2}W\left\{(d_{(1)}d_{(2)} +d_{(2)}d_{(1)}){\Re ic}\right\}\, .
$$
Now, it is very easy to produce an index version of this equation.  
(See \cite{AE1} for different versions). 

\

\section{Summary and Discussion}
The inspiration for this paper was the result --- that there exists a
single potential for any Weyl candidate in any dimension --- which we
obtained in \cite{ES} in a rather pragmatic manner; here, that result
has been shown to be a special case of a much more general result,
which itself is a consequence of a significant generalisation and
innovation in formalism.  The underlying approach has been the
systematic consideration of tensors as $r$-fold forms which has been
explained and discussed at length in \cite {S}.

The generalisation of the differential form approach to tensor-valued
forms and the extension of the use of the exterior differential $d$ to such
quantities is well known; we have emphasised the extension of the use
of $\delta$ and $\Delta$ also, and highlighted how an $r$-fold form
can be thought of as $r$ different tensor-valued forms, by taking each
of the $r$ sets of antisymmetric indices as the form indices, and
defining the three operators $d_{(i)}$, $\delta_{(i)}$ and $\Delta_{(i)},
\ i =1, \ldots r$, associated with each in turn.  However, we believe
the explicit introduction of the {\it weighted de Rham operator}
$\bar\D$ as defined in this paper adds an important new ingredient. 
Two crucial properties of this operator are that it enables us, via a
simple Laplace-like equation, to define an {\it associated
superpotential} of exactly the same tensor type, and in addition to
define {\it potentials in a very natural manner}; other generalised
Laplacian operators for tensor-valued forms lack one or both of these
properties.  Hence we believe that $\bar \D$ is a more powerful
operator than $\D_L$, and so will be very useful in the type of formal
investigations where $\D_L$ has been used previously \cite{Lic0},
\cite{Lic}, \cite{Bes}, and as a powerful alternative in
investigations of harmonic tensors. 

As a matter of fact, we have shown 
that the harmonicity property derived from the use of $\D_{L}$ is 
different from that derived from $\bar\D$, and that the latter has 
better properties concerning the harmonicity of Riemann curvature 
tensors (Section \ref{sec:Riemann}) and in proper Riemannian 
manifolds (Theorem \ref{harmonic=full}). We have also proved that one 
can use the defining properties of $\bar\D$ to define potentials and 
thereby to construct a generalised Hodge decomposition.
Theorem \ref{3} provides this fundamental generalisation to arbitrary tensors
(viewed as $r$-fold forms) of the local Hodge decomposition for
antisymmetric tensors (single $p$-forms); an obvious next question to
be considered is the possibility of a generalisation of the {\it
global} Hodge decomposition for closed proper Riemannian manifolds.  A
first partial step towards this direction has been established in
Theorem \ref{harmonic=full}, which we believe was not known hitherto. 

The remaining task will require, on the one hand, a systematic
analysis of the powers of the operators $d_{(i)}$ and $\delta_{(i)}$,
and on the other hand, a complete resolution of the relation between
tensors with the property $d_{(i)}T=0$, and those with the property
$\exists S:\,\, T=d_{(i)}S$ (and the dual analogue
with $\delta_{(i)}$ instead of $d_{(i)}$.) By the way, this will 
allow us to see if there are some implications of the Bianchi identities 
(\ref{bianchi},\ref{bianchi2}) for the potentials defined in 
Theorem \ref{theoremR}, in the same 
way as the first of the Maxwell equations (\ref{maxwell}) implies the existence 
of a  unique electromagnetic potential such that $F=dA$. All these matters 
are closely related,
actually, to the gauge problem for the general potentials that we have
mentioned in the paper and is the subject of current investigation
\cite{ESg}.  We are going to comment along these lines now.

As is well known, the operators $d$ and $ \delta$ are clearly
nilpotent for single $p$-forms in curved spaces, and for tensor-valued
forms in flat spaces.  Deeper nilpotent properties for tensor-valued
forms in curved space, and deeper commutator relationships with $\bar
\Delta$, and other operators (such as the first order operator which
defines the potential for the Weyl tensor) also need to be explored.

As an application of the general result for all tensors in Theorem
\ref{3}, we considered the Riemann curvature tensor showing that it
can be written in terms of a pair of potentials, as given in Theorem
\ref{theoremR}.  (There had earlier been some suggestions that the
Riemann tensor could be written in terms of {\it one} potential via the Riemann-Lanczos equations \cite{BU}, but
subsequently this was shown not to be possible \cite{Erh}, \cite{Er} --- at least for dimensions $n\ge 3$;
the fact that such a single potential cannot exist is now better
understood from this new result.)  The usefulness of this pair of
potentials, and in particular, the implications which the Bianchi
equations impose on their relationship, need further study.

Recently there have been attempts in {\it three} dimensions to obtain
$(2,1)$-form potentials $L_{ab}{}^c$ for the Ricci tensor via the {\it
Ricci-Lanczos equations} \cite{CG}, \bea R_{ab}= L_{aeb}{}^{;e} +
L_{bea}{}^{;e}+ L^e{}_{ae;b}+ L^e{}_{bea} \ .  \label{R=L} \eea The
reasons for the lack of success in obtaining such potentials for
arbitrary Ricci tensors in \cite{CG} become apparent when (\ref{R=L})
is compared with (\ref{Ric=Y+Z}); from the latter it is clear that
this type of potential formulation of the Ricci tensor requires {\it
two} potentials.  On the other hand, we know from (\ref{Rict=Y}) that
the {\it traceless part} of the Ricci tensor can always be given in
terms of {\it one} potential, in {\em arbitrary} dimensions; 
hence the partial success in \cite{CG} in
investigating special spaces in three dimensions {\it with zero Ricci
scalar}.  Moreover, we can see directly that, by excluding the Ricci
scalar from the potential structure, we will avoid the global
obstructions involving the Ricci scalar in \cite{CG} in all dimensions
(and the local obstruction in \cite{Er} in three dimensions).  It is
also instructive to compare the complicated nature of the second order
equations for the various superpotentials in \cite{CG} with the
corresponding simple and natural counterpart (\ref{barDelta11})
involving the weighted de Rham operator $\bar \D$.

There are many interesting aspects of the new potential for the Weyl
tensor which require deeper analysis.
When the new potential for the Weyl tensor in all dimensions was
obtained in \cite{ES} we pointed out that, {\it in four dimensions},
this new $(2,3)$-form potential $P^{ab}{}_{cde}$ is identical to the
double dual of the classical Lanczos $(2,1)$-form potential
$H^{ab}{}_{c}$ \cite {La}, \cite{I}, \cite{EH3}, \cite{AE2}.  It is
now clear why attempts to find higher dimensional analogues of the
Lanczos potential gave negative conclusions \cite{EH1}, \cite {EH2},
\cite{E3}; we should have been looking for higher dimensional
analogues not of the Lanczos potential, but of its double dual.  
There
do not seem to be deeper implications for this new result in four
dimensions; investigations of the Lanczos potential in four dimensions
seem to be most efficiently dealt with in spinors \cite{I},
\cite{AE2}.  However there are very important implications in higher
dimensions.

As emphasised in \cite{ES}, we now have, in all dimensions, an explicit
potential for the Weyl tensor which supplies a tensor which is an
`integral' of the Weyl tensor at the level of the connection, and
whose `square' has units $L^{-2}$ which are precisely the units we
would expect for gravitational `energies'. This suggests the
usefulness of the wave equation and the super-energy tensor of this
potential, {\it in all dimensions with Lorentz signature,} as an
alternative to the Bel-Robinson tensor for such mathematical
investigations as positivity properties, stability and the Cauchy
problem for the Einstein equations \cite{CK,AM}.

We have built our results on second order Laplace-like equations, with
appeals to the Cauchy-Kovalewski theorem \cite{CH} so that our results
are local for analytic pseudo-Riemannian metrics.  However, we expect
that these results can be generalised in the usual way; from the point
of view of general relativity we can appeal to stronger theorems
\cite{Fri}, \cite{HE} when we specialise to spaces with Lorentz
signature.  However, of course the potentials are by definition first
order, and now having established their existence it would be more
natural to consider them in a first order system.  In particular, the
single potential for the Weyl tensor is an attractive candidate for
deeper analysis in this context; preliminary investigations indicate
that more direct and powerful results can be obtained by treating this
definition as part of a first order symmetric hyperbolic system.

When we choose a potential we know that it is not unique, and a full
understanding of the role of gauge, which is more complicated for
tensor-valued forms than for single $p$-forms, will be very important
for further work.  In particular we will need to have a set of
explicit gauge equations to complete the first order symmetric
hyperbolic system for the Weyl tensor.  Furthermore, for the Weyl
tensor, we would hope for a second-order linear equation for its
potential with principal part of type $\nabla^h\nabla_h
P^{ab}{}_{cde}$, so that this will give an elliptic equation for
positive-definite metrics and a {wave equation} for Lorentzian
signature; in order to obtain such a simple version of the
second-order linear equation for the potential we will need to exploit
the gauge freedom.  In curved {\it four} dimensional space (where we
can be guided by the Lanczos potential \cite{I}, \cite{AE2},
\cite{El}), and in $n$-dimensional flat space it is easy to identify
this gauge freedom and obtain the expected simple Laplace-like
equation for the potential, but for other dimensions in curved
space the calculations become long and complicated \cite{ESg}.

\section*{Acknowledgements} 
JMMS gratefully acknowledges financial support from the Wenner-Gren Foundation, 
Sweden, and from grants FIS2004-01626 of the Spanish CICyT and 
no. 9/UPV 00172.310-14456/2002 of the University of the Basque 
Country. JMMS thanks the Matematiska institutionen, 
Link\"opings universitet, where this work was partly carried out, 
for hospitality.

\

\


\begin{thebibliography}{99}

\bibitem{AE1} Andersson, F. and Edgar, S.B. (1996).  The wave-equation for the Weyl tensor/spinor and dimensionally dependent tensor identities. {\it Int. J. of Mod. Physics}, {\bf 5}, 217-225.

 
\bibitem{AE2} Andersson, F. and Edgar, S.B. (2001).   Existence of Lanczos 
potentials and superpotentials for the Weyl spinor/tensor.  {\it 
Class. Quantum Gravity}, {\bf 18}, 2297-2304.
 
 
 
\bibitem{AM} Andersson, L. and Moncrief, V. (2003). 
Elliptic-hyperbolic systems and the Einstein equations. 
{\it Ann. Henri Poincar\'e} {\bf 4}, no. 1, 1--34.

%\bibitem{AE3} Andersson, F. and Edgar S.B. (2001).  Local existence of 
%symmetric spinor potentials for symmetric (3,1)-spinors in Einstein 
%space-times. {\it J. Geom. and Physics}, {\bf 37}, 273-290.

 
%\bibitem{AD}  Atkins, W. K. and Davis, W. R. (1980).  Lanczos-type potentials 
%in non-Abelian gauge theories. {\it Il Nuovo Cimento}, {\bf 59 B}, 116-132.

\bibitem{BF} Babourova, O.V. and Frolov, B.N.  (1995).  On a harmonic property of the Einstein manifold curvature. \  arXiv:gr-qc/9503045



\bibitem{Bel} Bel, L. (1963). \'Etude de certains op\'erateurs definis sur
les formes tensorielles (r,s). {\it  Ann. di Mat. Pura ed Applicata}, {\bf 51}, serie 
4,  171-192 


\bibitem{BT} Benn, I.M. and Tucker, R.W.  (1987).  {\it An Introduction to Spinors and Geometry with Applications in Physics}, Adam Hilger,  Institute of Physics.


\bibitem{Bes} Besse, A. L.   (1987).  {\it Einstein Manifolds}, Springer-Verlag.


%\bibitem{BC} Bampi, F., and Caviglia, G. (1983). Third-order tensor potentials 
%for the Riemann and Weyl tensors.  {\it Gen. Rel. Grav.}, {\bf 15}, 375.

\bibitem{BCJR}  Bini, D., Cherubini, C., Jantzen, R.T., and Rufini, R.  (2002).
{Teukolsky Master Equation: De Rham wave equation for the gravitational
and electromagnetic fields in vacuum.} {\it Prog.  Theor.  Phys.},
{\bf 107}, 967-992.


\bibitem{BCJF} Bini, D., Cherubini, C., Jantzen, R.T., and Rufini, R.   (2003).
{De Rham wave equations for tensor valued $p$-forms.} 
{\it Int. J. Mod. Phys. D}, {\bf 12}, 1363-1384.

\bibitem{BU} Brinis Udeschini, E. (1977). Su un tensore triplo potenziale del tensore di Riemann.
 {\it Mech. Fis. Mat. Istituto Lombardo (rend. Sc.)}, {\bf 111}, 466. 
\ (1980) \  A natural scalar field in the Einstein Gravitational Theory.
{\it Gen. Rel. Grav.}, {\bf 12}, 429.
 

\bibitem{Car} Cartan, \'E. (1946).  
{\it Lecons sur la g\'eom\'etrie des espaces de Riemann}, 
2nd edn. Gauthier-Villars, Paris.


\bibitem{CK} Christodoulou, D. and Klainerman, S. (1993).  
{\it The global nonlinear stability of Minkowski space}, 
Princeton Univ. Press.


\bibitem{CG} Chr\'usciel  P.T. and Gerber, A. (2005).  Some potentials for the curvature tensor on three-dimensional manifolds.  {\it Gen. Rel. Grav.}, {\bf 37},  891-905.



\bibitem{Coq} Coquereaux, R. (2002).  {\it Espaces, Fibr\'es et
Connexions},
http://www.cpt.univ-mrs.fr/\~{}coque/book/sourceforhtml.html




\bibitem{CH} Courant, R. and Hilbert, D. (1989).  {\it Methods of
Mathematical Physics Vol II}, Willey Classics Edition.


\bibitem{DR} de Rham, G. (1955).  {\it Vari\'et\'es
diff\'erentiables}, Hermann, Paris.

\bibitem{De} Deser, S. (1967). 
Covariant decomposition of symmetric tensors and the gravitational Cauchy problem. 
{\it Ann. Henri Poincar\'e} Section A, {\bf VII}, no. 2, 149--188.



\bibitem{DW} de Witt, B. (1962) in {\it Gravitation: An Introduction to Current Research}, ed. L. Witten, Witney, New York.


\bibitem{DH}  Dubois-Violette, M. and  Henneaux, M. (2002).
{Tensor fields of mixed Young symmetry type and 
N-complexes.} {\it Commun.Math.Phys}. {\bf 226}, 393-418.

\bibitem{Erh} Edgar, S. B.   (1994).  Non-existence of the Lanczos potential for 
the Riemann tensor in higher dimensions. 
 {\it Gen. Rel. Grav.}, {\bf 26},  329-332.



\bibitem{El} Edgar, S. B.  (1994).   The wave equation for the Lanczos 
tensor/spinor and a new tensor identity.  {\it Mod. Phys. Lett. A}, 
{\bf 9}, 479-482.



\bibitem{Er} Edgar, S. B. (2003). On effective constraints for the Riemann-Lanczos system of equations.
 {\it  J. Math. Phys.}, {\bf 44}, 5375-5385. 

 

\bibitem{E3}  Edgar, S. B. (2005).  Proofs of existence of local potentials for 
traceless symmetric 2-forms  using dimensionally dependent 
identities.   {\it J. Geom. Phys.}, {\bf 54}, 251-261. 
%{\it Preprint}: arXiv:gr-qc/0407002.

\bibitem{EH3} Edgar, S. B. and H\"oglund, A. (1997).  The Lanczos potential for 
the Weyl curvature tensor: existence, wave equations and algorithms.  
{\it Proc. Roy. Soc. A}, {\bf 453}, 835-851.


\bibitem{EH1} Edgar, S. B. and H\"oglund A.  (2000).  The Lanczos potential for 
Weyl candidates exists only in four dimensions. 
 {\it Gen. Rel. Grav.}, {\bf 32},  2307-2318.
 

%\bibitem{EH} Edgar, S.B. and H\"oglund, A.  (2002). Dimensionally dependent  tensor identities by %double antisymmetrisation.  {\it  J. Math. Phys.}, {\bf 43}, 659-677.


 



\bibitem{EH2}  Edgar, S. B. and H\"oglund, A. (2002).  The non-existence of a 
Lanczos potential for a Weyl curvature tensor in  dimensions $n\ge 
7$. 
 {\it Gen. Rel. Grav.}, {\bf 34}, 2149--2153.

%\bibitem{EW} Edgar, S.B. and Wingbrant, O.  (2003).  Old and new results for 
%superenergy tensors using dimensionally
%dependent identities.  {\it  J. Math. Phys.}, {\bf 44}, 6140-6159.


\bibitem{ES} Edgar, S. B. and Senovilla, J.M.M. (2004). A local potential 
for the Weyl tensor in all dimensions. {\it Class. Quantum Grav.}, {\bf  21},  
L133-L137.

\bibitem{ESg} Edgar, S. B. and Senovilla, J.M.M. (2005).  Gauge
freedom for general potentials and new potentials for curvature
tensors; in preparation.

\bibitem{Fl} Flanders, H. (1989). {\it Differential Forms with Applications to the Physical Sciences}, Dover Publications. 

\bibitem{Fr} Frankel, T. (1997). {\it The Geometry of Physics}, Cambridge University Press. 


\bibitem{Fri} Friedlander, F. G. (1975). {\it The Wave Equation on a Curved 
Spacetime}, Cambridge University Press.


\bibitem{Gem}  Gemelli, G. (2000).
{Second-order covariant tensor decomposition in curved spacetime.} {\it Mathematical Physics, Analysis and Geometry}, {\bf 3}, 195-216.

\bibitem{GS}  G\"ockeler, M. and Sch\"ucker, T. (1987).
{\it Differential geometry, gauge theories, and gravity.} Cambridge Monographs on Mathematical Physics, Cambridge University Press.

\bibitem{Gol} Goldberg, S.I. (1998 revised ed.) {\it Curvature and Homology},  Dover, London

\bibitem{HE} Hawking, S. W. and Ellis, G. F. R. 
(1973).   {\it The large scale structure of space-time},  
Cambridge Monographs on Mathematical Physics, No. 1. 
Cambridge University Press, London-New York.

\bibitem{H}  Helmholtz, H. (1858).
{\"Uber Integrale der hydrodynamischen Gleichungen, welche den
Wirbelbewegungen entsprechen.} {\it J. Reine Angew.  Math.}, {\bf 55},
25-55

\bibitem{Ho} Hodge,  W. V. D. (1941) {\it The Theory and Applications of Harmonic 
Integrals}, Cambridge University Press, Cambridge.

\bibitem{I} Illge, R. (1988).  On potentials for several classes of spinor and 
tensor fields in curved space-times.  {\it Gen. Rel. Grav.}, {\bf 20}, 551-564.

\bibitem{Kul} Kulkarni, R. S. (1972).  On the Bianchi Identities.  {\it Math. Ann.}, {\bf 199}, 175-204.

\bibitem{La} Lanczos, C. (1962). The splitting of the Riemann tensor. {\it  
Rev. Mod. Phys.}, {\bf 34}, 379-389.

\bibitem{Lic0} Lichnerowicz, A. (1961) {\it Propagateurs et 
commutateurs en relativit\'e g\'en\'erale}, Publ. Scient. des Hautes 
\'etudes scientifiques, No.10, Hermann, Paris.

\bibitem{Lic} Lichnerowicz, A.  (1964).  {\it Propagateurs, commutateurs, et 
anticommutateurs en relativit\'e g\'en\'erale}, 
in ``Relativity, Groups and Topology'', eds. C. DeWitt and B.S. 
DeWitt, Gordon and Breach, New York.

\bibitem{LD} Lindell, L. V. and Dassios, G. (2003).   Helmoltz Theorem for multiform fields.  {\it 
J. Electromagn. Waves and Appl.,} {\bf 17}, 3-14.


%\bibitem{Lov} Lovelock, D. (1970).   Dimensionally dependent identities.  {\it 
%Proc. Camb. Phil. Soc.}, {\bf 68}, 345-350.


\bibitem{MTW} Misner, C.W., Thorne, K. S. and Wheeler, J.A. (1973). {\it Gravitation}, Freeman.

%\bibitem{R} Roberts, M. D. (1989). Dimensional reduction and the Lanczos 
%tensor.  {\it Mod. Phys. Letts. A}, {\bf 4}, 2739-2746.

\bibitem{Na} Nakahara, M. (1990). {\it Geometry, Topology and Physics}, Institute of Rhysics.


\bibitem{MT} Parker, L. and Christensen, S.M.  (1994). {\it MathTensor},  Addison-Wesley, New York. 

\bibitem{Pe} Penrose, R. (1960). {\it Ann. Phys.}, {\bf 10}, 171-201.


\bibitem{Ry} Ryan, M.P. (1974). Teukolsky equation and Penrose wave equation. {\it Phys. Rev., D}, {\bf 10}, 1736-1740.

\bibitem{S} Senovilla, J. M. M.  (2000).  Super-energy Tensors.   {\it 
Class. Quantum Grav.}, {\bf 17}, 2799-2841.

\bibitem{Sc} Schwarz, G. (1991). {\it Hodge Decomposition --- a Method for Solving Boundary value Problems}, Springer-Verlag.


\bibitem{Th} Thirring, W. (1978). {\it Classical Field Theory}, Springer-Verlag.

\bibitem{Ur} Urakawa, H. (1993). {\it Calculus of Variations and Harmonic Maps}, A.M.S..


\bibitem{W} Woodside, D. A.  (1999). Uniqueness theorems for classical four-vector fields in Euclidean and Minkowski spaces.  {\it  J. Math. Phys.}, {\bf 
40}, 4911-4943.

\bibitem{Yo} York, J.W. (1973). Conformally invariant orthogonal decomposition of symmetric tensors on Riemannian manifolds and the initial-value problem of general relativity. {\it J. Math. Phys.}, {\bf 14}, 456-471.  

\end{thebibliography}
\end{document}